\documentclass[12pt,a4paper]{article}
\usepackage{amstext}
\usepackage{fancyhdr}
\usepackage{amsfonts,graphicx,bezier, amssymb}
\usepackage{amsmath}
\usepackage{caption}
\usepackage{mathtools}
\usepackage{mathrsfs}
\usepackage[utf8]{inputenc}
\usepackage{tikz}
\usepackage{tikz-cd}
\usepackage{rotating}
\usepackage[multiple]{footmisc}
\usetikzlibrary{arrows.meta}
\newenvironment{prf}{\noindent{\bf{Proof:}}~~}{\hfill\rule{1ex}{1ex}\vskip1.5ex}
\newcommand{\Z}{\mathbb Z}

\newcommand{\beqa}{\begin{eqnarray}}
\newcommand{\enqa}{\end{eqnarray}}
\newcommand{\beq}{\begin{eqnarray*}}
	\newcommand{\enq}{\end{eqnarray*}}

\newtheorem{rem}{Remark}[section]

\newtheorem{cor}{Corollary}[section]
\newtheorem{prop}{Proposition}[section]
\newtheorem{defn}{Definition}[section]
\newtheorem{exam}{{\bf Example}}[section]
\newtheorem{thm}{Theorem}[section]

\newtheorem{lem}{Lemma}[section]

\newcommand{\noi}{\noindent}

\makeatletter
\providecommand*{\twoheadrightarrowfill@}{%
	\arrowfill@\relbar\relbar\twoheadrightarrow
}
\providecommand*{\twoheadleftarrowfill@}{%
	\arrowfill@\twoheadleftarrow\relbar\relbar
}
\providecommand*{\xtwoheadrightarrow}[2][]{%
	\ext@arrow 0579\twoheadrightarrowfill@{#1}{#2}%
}
\providecommand*{\xtwoheadleftarrow}[2][]{%
	\ext@arrow 5097\twoheadleftarrowfill@{#1}{#2}%
}
\makeatother
\usepackage[para]{manyfoot}
\begin{document}
	\begin{center}
		{\bf\Large On the Greenlees-May Duality and the Matlis-Greenlees-May Equivalence}
		%Categorical Generalization of Reduced and Coreduced modules
	\end{center}
		\vspace*{0.3cm}
		\begin{center}
		Tilahun Abebaw,\footnote{Department of Mathematics, Addis Ababa University, P.O.Box 1176, Addis Ababa, Ethiopia \newline email: tilahun.abebaw@aau.edu.et, amanuel.mamo@aau.edu.et, zelalem.teshome@aau.edu.et}  Amanuel Mamo,\footnotemark[1]{}\footnote{This work forms part of the second author's PhD thesis.}
	    %\DeclareNewFootnote{A}[]
	    David Ssevviiri  \footnote {Department of Mathematics, Makerere University, P.O.Box 7062, Kampala Uganda, email: david.ssevviiri@mak.ac.ug}\footnote{Corresponding author} and
	      Zelalem Teshome\footnotemark[1]	
			\vspace*{0.3cm}
		\end{center}
		\begin{abstract}\noi 
		Let $A$ be a commutative unital ring and $\mathfrak{ a}$ an ideal in it. We define and study $\mathfrak{ a}$-reduced complexes and $\mathfrak{ a}$-coreduced complexes in both the category of chain complexes of $A$-modules as well as in the derived category of $A$-modules. We show that these two types of complexes give rise to variants of the well known Greenlees-May duality and the Matlis-Greenlees-May equivalence in the aforementioned categories. 
		\end{abstract}
		\vspace*{0.3cm}

		{\bf MSC 2020} Mathematics Subject Classification: 13D07, 13D09, 18A40, 18G35,  18A35. 
		\section{Introduction}
			\begin{paragraph}\noi 
			Let $A$ be a commutative unital ring, and $\mathfrak{a}$ an ideal of $A$. Consider the following endo functors on the category of $A$-modules; the $\mathfrak{a}$-torsion functor $$\Gamma_{\mathfrak{a}}(-):=\underset{k}{\varinjlim}\text{Hom}_{A}(A/\mathfrak{a}^{k},-)$$ and the $\mathfrak{a}$-adic completion functor $$\Lambda_{\mathfrak{a}}(-):=\underset{k}{\varprojlim}(A/\mathfrak{a}^{k}\otimes_{A}-).$$ These functors are dual to each other and have been utilised, among other reasons, to define $\mathfrak{a}$-torsion modules and $\mathfrak{a}$-adically complete modules respectively. An $A$-module $M$ is \textit{ $\mathfrak{ a}$-torsion} (resp.\textit{ $\mathfrak{ a}$-adically complete}) if $\Gamma_{ \mathfrak{a}}(M)\cong M$ (resp. $M\cong \Lambda_{ \mathfrak{a}}(M)$). $\mathfrak{a}$-torsion modules were characterised in terms of functors by Positselski in \cite[Sections 6 and 7]{Leonid_contraadjusted...}, and $\mathfrak{a}$-adically complete modules, on the other hand, were characterised in terms of functors in \cite[Chapter $3$]{Peter-Schenzez}.
%			 For $s\in A$, $s$-(adically) complete modules form a subclass of $s$-contramodules studied in \cite[Theorem 2.3]{Leonid_contraadjusted...}. 
 The homological properties of these two classes of modules were studied in \cite{porta2014homology, Peter-Schenzez}, among others. This was achieved by defining and studying  derived $\mathfrak{a}$-torsion complexes and derived $\mathfrak{a}$-adically complete complexes. These complexes are alternatively referred to as cohomologically $\mathfrak{a}$-torsion complexes and (co)homologically $\mathfrak{a}$-(adically) complete complexes in the literature respectively, see for instance, \cite[Section 1]{porta2014homology} and \cite[Definition 9.6.1]{Peter-Schenzez}.   The category of derived $\mathfrak{a}$-torsion complexes is equivalent to the category of derived $\mathfrak{a}$-adically complete complexes, provided $\mathfrak{ a}$ is a weakly proregular ideal, see \cite[Theorem 7.11]{porta2014homology} and \cite[Theorem 9.6.7]{Peter-Schenzez}. This is what is called the \textit {Matlis-Greenlees-May equivalence} (MGM equivalence for short). The equivalence is established by utilising the right derived functor $\text{R}\Gamma_{ \mathfrak{a}}$ of $\Gamma_{ \mathfrak{a}}$, and the left derived functor $\text{L}\Lambda_{ \mathfrak{a}}$ of $\Lambda_{ \mathfrak{a}}$. A version of the equivalence for complexes over a non-commutative ring was given in \cite{vyas2018weak}.  Moreover,
		 for ideals $\mathfrak{ a}$ and $\mathfrak{ b}$ of a commutative Noetherian ring $A$, the MGM equivalence is known to exist between the category  of derived $(\mathfrak{a},\mathfrak{b})$-torsion complexes and the category of derived $(\mathfrak{a},\mathfrak{b})$-adically complete complexes \cite{Li-Yang-Generalised MGM}. 
 \end{paragraph}
		\begin{paragraph}\noi
			The other important property about the derived functors $\text{R}\Gamma_{\mathfrak{a}}$ and $\text{L}\Lambda_{\mathfrak{a}}$ was demonstrated in \cite[Theorem 7.12]{porta2014homology} and \cite[Corollary 9.2.4]{Peter-Schenzez}, among others. It states that for a weakly proregular ideal $\mathfrak{a}$ of $A$, the functor $\text{R}\Gamma_{\mathfrak{a}}$ is derived left adjoint to $\text{L}\Lambda_{\mathfrak{a}}$, i.e., $\text{RHom}_{A}(\text{R}\Gamma_{ \mathfrak{a}}(M),N)\cong \text{RHom}_{A}(M, \text{L}\Lambda_{ \mathfrak{a}}(N))$ for any complexes $M$ and $N$. This is what is called the \textit{Greenlees-May duality} (GM duality for short).
			The GM duality is a far-reaching generalisation of the celebrated Grothendieck local duality \cite{Alonso-Lipman}. Moreover, both the MGM equivalence and the GM duality were first established in \cite{Alonso-Lipman} in the algebraic geometry setting of schemes.
			
		\end{paragraph} 
		\begin{paragraph}\noi In this paper, we utilise the functors $\Gamma_{\mathfrak{a}}$ and $\Lambda_{\mathfrak{a}}$ (resp. their derived functors $\text{R}\Gamma_{ \mathfrak{a}}$ and $\text{L}\Lambda_{ \mathfrak{a}}$) to define $\mathfrak{a}$-reduced complexes and $\mathfrak{a}$-coreduced complexes in the category of chain complexes of $A$-modules, $\text{C}(A)$, (resp. in the derived category of $A$-modules, $\text{D}(A)$).  We characterise $\mathfrak{ a}$-reduced complexes and $\mathfrak{ a}$-coreduced complexes, and examine their homological properties. In particular, we provide versions of the GM duality and the MGM equivalence in both $\text{C}(A)$ and $\text{D}(A)$ using $\mathfrak{a}$-reduced complexes and $\mathfrak{a}$-coreduced complexes and without assuming that $\mathfrak{ a}$ is a weakly proregular ideal of $A$. To be specific,  our main results are Theorems \ref{Thm 1.1}, \ref{GMD in C(A)}, \ref{Prop: version of GM Duality} and  \ref{Thm: MGM Equality} below which appear in the body of the paper as Theorems \ref{sec3: GM Duality-Thm}, \ref{MGM in C(A)}, \ref{Thm-version of GM Duality in D(A)} and  \ref{Theorem- MGM equality} respectively. It should be noted that this is an extension of results proved in \cite[Theorems 3.4 and 4.3]{David-Application-I} in the category of $A$-modules.
		
		\end{paragraph}  
		
		\begin{paragraph}\noi To understand the four theorems, we make the notion of $\mathfrak{ a}$-reduced and $\mathfrak{ a}$-coreduced complexes precise first and elaborate a bit how they evolved from the notion of reduced rings.
%		\end{paragraph}
%\begin{paragraph}\noi
	A ring is \textit{reduced} if it has no nonzero nilpotent elements. This notion was extended to modules by Lee and Zhou in \cite{Lee and Zhou}. An $A$-module $M$ is \textit{reduced} if for all ideals $\mathfrak{ a}$ of $A$ and all $m\in M$, $\mathfrak{ a}^{2}m=0$ implies that $\mathfrak{ a}m=0$. We ``localise'' this concept by defining $\mathfrak{ a}$-reduced modules for a fixed ideal $\mathfrak{ a}$ of $A$. For an ideal $\mathfrak{ a}$ of $A$, an $A$-module $M$ is \textit{$\mathfrak{ a}$-reduced} if for all $m\in M$, the implication $\mathfrak{ a}^{2}m=0\Rightarrow \mathfrak{ a}m=0$ holds. It was shown in \cite[Proposition 2.2]{David-Application-I} that $M$ is $\mathfrak{ a}$-reduced if and only if $\mathfrak{ a}\Gamma_{ \mathfrak{a}}(M)=0$. Dually, an $A$-module $M$ is \textit{$\mathfrak{ a}$-coreduced} if and only if $\mathfrak{ a}M=\mathfrak{ a}^{2}M$. This is also equivalent to saying that $\mathfrak{ a}\Lambda_{ \mathfrak{a}}(M)=0$ \cite[Proposition 2.3]{David-Application-I}.  A complex of $A$-modules is  \textit{$\mathfrak{ a}$-reduced} (resp. \textit{$\mathfrak{ a}$-coreduced}) if each module in it is $\mathfrak{ a}$-reduced (resp. $\mathfrak{ a}$-coreduced). It is worth noting that $\mathfrak{ a}$-reduced $A$-modules, $\mathfrak{ a}$-coreduced $A$-modules and their respective generalisations have been recently utilised for several purposes. For details, see \cite{Tekle-David, reduced w.r.t another-ours, Kimuli-David,  kyomuhangi2020locally,Annet-David : Generalized reduced, David-Application-I, ssevviiri: App II}. 
	%We prove that $\mathfrak{ a}$-reduced complexes and $\mathfrak{a}$-coreduced complexes provide convenient conditions for realisation of versions of GM duality and the MGM equivalence. In particular, we prove the following four main theorems. 
	\end{paragraph}
		\begin{paragraph}\noi
			For an ideal $\mathfrak{ a}$ of a ring $A$, denote by $\text{C}(A)_{\mathfrak{a}\text{-red}}$ and  $\text{C}(A)_{\mathfrak{a}\text{-cor}}$ the full subcategory of $\text{C}(A)$ consisting of $\mathfrak{a}$-reduced complexes and $\mathfrak{a}$-coreduced complexes respectively.
		\end{paragraph}
		\begin{thm}[The version of GM duality in $\text{C}(A)$]\label{Thm 1.1}
			Let $A$ be a ring and $\mathfrak{a}$ an ideal in it.
			\begin{itemize}
				\item[$(1)$]  The functor
				{\normalfont\begin{equation*}
						\Gamma_{\mathfrak{a}}: \text{C}(A)_{\mathfrak{a}\text{-red}}\rightarrow \text{C}(A)_{\mathfrak{a}\text{-cor}}; ~ M\mapsto \Gamma_{\mathfrak{a}}(M)
				\end{equation*}}
				is idempotent and {\normalfont $\Gamma_{\mathfrak{a}}(M)\cong \text{Hom}_{A}(A/\mathfrak{a},M)$} for any $\mathfrak{a}$-reduced complex $M$.
				\item[$(2)$] The functor 
				{\normalfont \begin{equation*}
						\Lambda_{\mathfrak{a}}: \text{C}(A)_{\mathfrak{a}\text{-cor}}\rightarrow \text{C}(A)_{\mathfrak{a}\text{-red}};~ M\mapsto \Lambda_{\mathfrak{a}}(M)
				\end{equation*}}
				is idempotent and {\normalfont $\Lambda_{\mathfrak{a}}(M)\cong A/\mathfrak{a}\otimes_{A} M$} for any $\mathfrak{a}$-coreduced complex $M$.
				\item[$(3)$]  For any $M\in {\normalfont \text{C}(A)_{\mathfrak{a}\text{-cor}}} $ and any {\normalfont $ N\in \text{C}(A)_{\mathfrak{a}\text{-red}}$}, 
				{\normalfont \begin{equation*}
						\text{Hom}_{A}(\Lambda_{\mathfrak{a}}(M) ,N)\cong 	\text{Hom}_{A}( M,\Gamma_{\mathfrak{a}}(N)).
				\end{equation*}}
			\end{itemize}
		\end{thm}
	
	\begin{paragraph}\noi  We call a complex \textit{ $\mathfrak{ a}$-torsion} (resp. \textit{$\mathfrak{ a}$-adically complete}) if each of its modules is $\mathfrak{ a}$-torsion (resp. $\mathfrak{ a}$-adically complete). Denote by $\text{C}(A)_{\mathfrak{a}\text{-tor}}$ and $\text{C}(A)_{\mathfrak{a}\text{-com}}$ the full subcategory  of $\text{C}(A)$ consisting of $\mathfrak{a}$-torsion complexes and $\mathfrak{a}$-adically complete complexes respectively.
		\end{paragraph}

		\begin{thm}[The version of MGM equivalence in $\text{C}(A)$]\label{GMD in C(A)} For any ideal $\mathfrak{ a}$ of $A$, we have the equalities:
			{\normalfont $$\text{C}(A)_{\mathfrak{a}\text{-tor}}\cap \text{C}(A)_{\mathfrak{a}\text{-red}}= \text{C}(A)_{\mathfrak{a}\text{-com}}\cap \text{C}(A)_{\mathfrak{a}\text{-cor}}=\{M\in\text{C}(A): \mathfrak{ a}M\cong 0\}.$$}
		\end{thm}
		\begin{paragraph}\noi 
			 A complex $M\in\text{D}^{}(A)$ is \textit {$\mathfrak{ a}$-reduced} (resp. \textit {$\mathfrak{ a}$-coreduced}) if $\text{RHom}_{A}(A/\mathfrak{ a}, M)\cong\text{R}\Gamma_{ \mathfrak{a}}(M)$ (resp. $\text{L}\Lambda_{ \mathfrak{a}}(M)\cong A/\mathfrak{ a}\otimes_{A}^{\text{L}}M$). We denote by $\text{D}(A)_{\mathfrak{ a\text{-red}}}$ (resp. $\text{D}(A)_{\mathfrak{ a\text{-cor}}}$) the  $\mathfrak{ a}$-reduced  (resp. $\mathfrak{ a}$-coreduced) complexes in $\text{D}(A)$.
		\end{paragraph}
		\begin{thm}[The version of GM duality in $\text{D}(A)$]\label{Prop: version of GM Duality} Let $\mathfrak{ a}$ be an ideal of a ring $A$, {\normalfont $M\in\text{D}^{}(A)_{\mathfrak{ a\text{-cor}}}$} and {\normalfont $N\in\text{D}^{}(A)_{\mathfrak{ a\text{-red}}}$}. There exists an isomorphism in {\normalfont $\text{D}(A)$} given by  {\normalfont $$ \text{RHom}_{A}(\text{L}\Lambda_{\mathfrak{ a}}(M), N)\cong \text{RHom}_{A}(M,\text{R}\Gamma_{\mathfrak{ a}}(N)).$$}	
		\end{thm}
		Let $\text{D}_\text{f}^\text{b}(A)$ denote the bounded complexes of finitely generated $A$-modules in $\text{D}(A)$, and $\text{D}_\text{f}^\text{b}(A)_{\mathfrak{ a\text{-red}}}$ the category of bounded complexes of finitely generated $\mathfrak{ a}$-reduced modules. Analogously, $\text{D}_\text{f}^\text{b}(A)_{\star}$ is defined, where $\star= \mathfrak{ a\text{-cor}},\mathfrak{ a\text{-tor}}, \mathfrak{ a\text{-com}}$.
		\begin{thm}[The version of MGM equivalence in $\text{D}(A)$] \label{Thm: MGM Equality}
			Let $A$ be a Noetherian ring such that each complex in {\normalfont $\text{D}_\text{f}^\text{b}(A)$} has both finite injective dimension and finite Tor dimension. We have the equality:
			{\normalfont \begin{equation*}
					\text{D}_\text{f}^\text{b}(A)_{\mathfrak{a}\text{-com}}\cap \text{D}_\text{f}^\text{b}(A)_{\mathfrak{a}\text{-cor}}= \text{D}_\text{f}^\text{b}(A)_{\mathfrak{a}\text{-tor}}\cap \text{D}_\text{f}^\text{b}(A)_{\mathfrak{a}\text{-red}}.
			\end{equation*}}
		\end{thm}
	
	\begin{paragraph}\noi 
		The paper is organised as follows.  In Section 2, we characterise and give several properties of $\mathfrak{a}$-reduced complexes and $\mathfrak{a}$-coreduced complexes.  Section 3 is devoted to proving versions of the GM duality and the MGM equivalence in $\text{C}(A)$. As a consequence, we show that if the ideal $\mathfrak{ a}$ is idempotent, then $\text{C}(A)_{\mathfrak{ a\text{-com}}}$, which is not abelian in general, turns out to be an abelian subcategory of $\text{C}(A)$. Finally, in Section 4, we characterise $\mathfrak{ a}$-reduced and $\mathfrak{ a}$-coreduced complexes and give versions of the GM duality and MGM equivalence in $\text{D}(A)$.
	\end{paragraph}
\section{Preliminaries and basic properties}
\begin{paragraph}\noi
Throughout the paper $A$ is a commutative ring with unity, and $\mathfrak{ a}$ is an ideal in $A$. In this section, we introduce $\mathfrak{a}$-reduced and $\mathfrak{a}$-coreduced  complexes over $A$ and study their basic properties. We denote by $\text{M}(A)$ the category of $A$-modules, by $\text{C}(A)$ its category of chain complexes and by $\text{D}(A)$ the derived category of $\text{M}(A)$.
\end{paragraph}
\begin{paragraph}\noi
We begin by recalling some definitions and facts to be used in the sequel. 
%Let $A$ be a ring and let $\mathfrak{ a}$ be an ideal of $A$.	
The complexes $ M:=(M^{i},d_{M}^{i})_{i\in\mathbb{Z}}$ and $ N:=(N^{j},d_{N}^{j})_{j\in\mathbb{Z}} \in \text{C}(A)$ are isomorphic if all the homomorphisms $f^{i}: M^{i}\rightarrow N^{i}$ are isomorphisms. When $M$ and $N$ are isomorphic we write $M\cong N$. 
Following \cite{Peter-Schenzez, Yekutieli-Derived categories} the homomorphism complex $\text{Hom}_{A}(M,N)$ and the tensor product complex $M\otimes_{A}N$ for complexes $ M=(M^{i},d_{M}^{i})_{i\in\mathbb{Z}}$ and $ N=(N^{j},d_{N}^{j})_{j\in\mathbb{Z}} \in \text{C}(A)$ are respectively given by
	$$ \text{Hom}_{A}(M,N):=(\text{Hom}_{A}(M, N )^ {n}, d^{n})_{n\in\mathbb{Z}},$$~ \text{where}	$\text{Hom}_{A}(M, N )^ {n}=\underset{j-i=n}{\prod} {\text{Hom}}_{A} (M^{i} , N^ {j})$
	and
	$d^{n}= \underset{j-i=n}{\prod}{\text{Hom}}_{A}(d_{M}^ {i-1},N^{j}) + (-1)^{i+1}\text{Hom}_{A}(M^{i},d_{N}^{j})$ and
	$$M\otimes_{A} N:= ((M\otimes_{A} N)^{n},\sigma^{n})_{n\in\mathbb{Z}},$$ ~\text{where}
	$(M\otimes_{A}N)^{n}=\underset{i+j=n}{\oplus}M^{i}\otimes_{A} N^{j}$
	and 
	$\sigma^{n}:= d^{i}_{M}\otimes_{A} \text{id}_{N} +(-1)^{i}~ \text{id}_{M}\otimes_{A} d^{i}_{N}$.		
\end{paragraph}
\begin{lem}\label{sec2: Hom-tensor Lem} Let $\mathfrak{a}$ be an ideal of $A$, and {\normalfont $M\in\text{C}(A)$}. For any $n\in\mathbb{Z}$,
	{\normalfont	\begin{itemize}
			\item [$(1)$]  $\text{ Hom}_{A}(A/\mathfrak{ a}, M)^{n}=\text{ Hom}_{A}(A/\mathfrak{ a}, M^{n})$.
			\item[$(2)$]	$(A/\mathfrak{ a}\otimes_{A} M)^{n}= A/\mathfrak{ a}\otimes_{A} M^{n}$.	
	\end{itemize}}
\end{lem}
\begin{prf}
	Let $M:=(M^{n})_{n\in\mathbb{Z}}$. Since for any $k\in\Z^{+}$ the $A$-module $(A/\mathfrak{ a})^{k}$ is an $A$-complex concentrated in degree zero, we get
	\begin{itemize}
	\item [$(1)$] {\normalfont	\begin{equation*}
		\text{Hom}_{A}(A/\mathfrak{a},M)^{n}= \underset{j-i=n}{\prod}\text{Hom}_{A}((A/\mathfrak{a})^{i},M^{j})=\underset{j=n}{\prod}\text{Hom}_{A}(A/\mathfrak{a},M^{j})
		\end{equation*}
		\begin{equation*}
		=\text{Hom}_{A}(A/\mathfrak{a},M^{n}).
		\end{equation*}}
	\item[$(2)$] 
	\begin{equation*}
	(A/\mathfrak{ a}\otimes_{A}M)^{n}=\underset{i+j=n}{\bigoplus} ((A/\mathfrak{ a})^{i}\otimes_{A} M^{j})= \underset{j=n}{\bigoplus} (A/\mathfrak{ a}\otimes_{A} M^{j})
	\end{equation*}
	\begin{equation*}
	=A/\mathfrak{ a}\otimes_{A} M^{n}.
	\end{equation*}
	\end{itemize}
\end{prf}
%\begin{defn}
%	{\normalfont Let $\mathfrak{a}$ be an ideal of $A$. An $A$-module $M$ is} $\mathfrak{a}$-torsion {\normalfont if $\Gamma_{\mathfrak{a}}(M)$ is isomorphic to $M$} {\normalfont and $M$ is called} $\mathfrak{ a}$-adically complete {\normalfont if $M$ is isomorphic to $\Lambda_{ \mathfrak{a}}(M)$, the $\mathfrak{ a}$-adic completion of $M$.}
%\end{defn}
%We simply denote by $M$ the complex $M:= (M^{n}, d_{M}^{n})_{n\in\mathbb{Z}}$.
\begin{defn}\label{Defn: reduced complex} {\normalfont A complex $ M$ is} $\mathfrak{a}$-reduced{ \normalfont  if each module in $M$ is an $\mathfrak{a}$-reduced $A$-module}.
\end{defn}
\begin{paragraph}\noi
%	A subcomplex of an $\mathfrak{a}$-reduced complex is $\mathfrak{a}$-reduced.	In particular, if $M$ is an $\mathfrak{a}$-reduced complex, then the subcomplex $\Gamma_{\mathfrak{a}}(M)$ is $\mathfrak{a}$-reduced. 
We denote by $\text{C}(A)_{\mathfrak{a}\text{-red}}$ the full subcategory of $\text{C}(A)$ consisting of  $\mathfrak{a}$-reduced complexes. 
\end{paragraph}
\begin{prop}\label{sec2: Pro-I-red} {\normalfont Let $\mathfrak{ a}$ be an ideal of $A$ and let $M$ be a complex over $A$. The following statements are equivalent:
		\begin{itemize} 
			\item[$(1)$] $ M\in \text{C}(A)_{\mathfrak{a}\text{-red}}$,  
			\item[$(2)$]$\text{Hom}_{A}(A/\mathfrak{a}, M)\cong \text{Hom}_{A}(A/\mathfrak{a}^{2},M)$,
			\item[$(3)$]$\Gamma_{\mathfrak{a}}(M)\cong$ $\text{Hom}_{A} (A/\mathfrak{a},M)$, 
			\item[$ (4)$] $\mathfrak{a}\Gamma_{\mathfrak{a}}(M)\cong 0$. 
	\end{itemize}}
\end{prop}
\begin{prf} 
	\begin{itemize}
		\item[$(1)\Rightarrow (2)$] Since $M=(M^{n})_{n\in\mathbb{Z}}$ is $\mathfrak{a}$-reduced for each $n$, and the complex $A/\mathfrak{a}^{k}, k\in\mathbb{Z^{+}}$ is concentrated in degree zero, by Lemma \ref{sec2: Hom-tensor Lem} $(1)$, we have that 
		$	\text{Hom}_{A}(A/\mathfrak{a}, M)^{n}=\text{Hom}_{A}(A/\mathfrak{a}, M^{n})
		\cong \text{Hom}_{A}(A/\mathfrak{a}^{2}, M^{n})
		= \text{Hom}_{A}(A/\mathfrak{a}^{2}, M)^{n}$. So, the complexes $ \text{Hom}_{A}(A/\mathfrak{a}, M) $ and $\text{Hom}_{A}(A/\mathfrak{a}^{2}, M)$ are isomorphic.
		\item[$ (2) \Rightarrow (3)$] From $(2)$, there is an isomorphism $\text{Hom}_{A}(A/\mathfrak{a},M)^{n}\cong \text{Hom}_{A}(A/\mathfrak{a}^{k}, M)^{n}$ of $A$-modules for any $k\in\mathbb{Z^{+}}$ and any
		$ n\in\mathbb{Z}$. By \cite[Proposition 2.2]{David-Application-I}, 
		$\Gamma_{\mathfrak{a}}(M)^{n}:=\underset{k}\varinjlim~ \text{Hom}_{A}(A/\mathfrak{a}^{k}, M)^{n}\cong \text{Hom}_{A}(A/\mathfrak{a},M)^{n}$  for every $n$.
		So, the complexes $\Gamma_{\mathfrak{a}}(M)$ and $\text{Hom}_{A}(A/\mathfrak{a},M)$ are isomorphic.
		\item[$ (3) \Rightarrow (4)$] Let $n\in\mathbb{Z}$. The  $A$-modules $\text{Hom}_{A}(A/\mathfrak{a},M)^{n}$ in the complex $\text{Hom}_{A}(A/\mathfrak{a},M)$ are annihilated by $\mathfrak{a}$. The conclusion is that $\mathfrak{a}\text{Hom}_{A}(A/\mathfrak{a},M)\cong \mathfrak{a}\Gamma_{\mathfrak{a}}(M)\cong 0$.
		\item[$(4)\Rightarrow (1)$]  Since $\mathfrak{a}\Gamma_{\mathfrak{a}}(M)\cong 0$, the $A$-modules $\mathfrak{a}\Gamma_{\mathfrak{a}}(M)^{n}\cong \mathfrak{ a}\Gamma_{ \mathfrak{a}}(M^{n})=0$ for all $n$. 
	%	Note that $\Gamma_{ \mathfrak{a}}(M^{n})$ is the $n^\text{th}$ component of the induced subcomplex $\Gamma_{ \mathfrak{a}}(M)$.
		It follows that $M^{n}$ is $\mathfrak{a}$-reduced for each $n$. Therefore, $M\in\text{C}(A)_{\mathfrak{ a\text{-red}}}$.
	\end{itemize}
\end{prf}
\begin{defn}\label{Defn: coreduced complexes}
	{\normalfont A complex is } $\mathfrak{a}$-coreduced {\normalfont if each of its $A$-module is $\mathfrak{a}$-coreduced.}
\end{defn}
Denote by $\text{C}(A)_{\mathfrak{a}\text{-cor}}$ the full subcategory of $\text{C}(A)$ consisting of $\mathfrak{a}$-coreduced complexes. We give  the dual of Proposition \ref{sec2: Pro-I-red}. 
\begin{prop}\label{sec2: Pro I-cor}
	Let $\mathfrak{a}$ be an ideal of $A$ and $M$ a complex of $A$-modules. The following statements are equivalent:
	\begin{itemize}
		\item [$(1)$] $M\in {\normalfont \text{C}(A)_{\mathfrak{a}\text{-cor}}}$,
		\item[$(2)$]$A/\mathfrak{a}\otimes_{A} M\cong A/\mathfrak{a}^{2}\otimes_{A} M$,
		\item[$(3)$]  $\Lambda_{\mathfrak{a}}(M)\cong A/\mathfrak{a}\otimes_{A} M$,
		\item[$(4)$] $\mathfrak{a}\Lambda_{\mathfrak{a}}(M)\cong 0$. 
	\end{itemize}
\end{prop}
\begin{prf}
	\begin{itemize}
		\item [$(1)\Rightarrow (2)$] Let $M=(M^{n})_{n\in\mathbb{Z}}$ be an $\mathfrak{a}$-coreduced complex. For each $n$, the $A$-modules $M^{n}$ are $\mathfrak{a}$-coreduced, and since $A/\mathfrak{a}^{k}, k\in\mathbb{Z^{+}}$ is concentrated in degree zero,
		$(A/\mathfrak{a}\otimes_{A}M)^{n}= A/\mathfrak{a}\otimes_{A} M^{n}\cong A/\mathfrak{a}^{2}\otimes_{A} M^{n}= (A/\mathfrak{a}^{2}\otimes_{A}M)^{n}$ [Lemma \ref{sec2: Hom-tensor Lem} $(2)$].		
		So, $(A/\mathfrak{a}\otimes_{A}M)^{n}\cong (A/\mathfrak{a}^{2}\otimes_{A}M)^{n}$. Therefore, the complexes $A/\mathfrak{a}\otimes_{A}M$ and $A/\mathfrak{a}^{2}\otimes_{A}M$ are isomorphic.
		\item[$(2)\Rightarrow (3)$] By $(2)$, for each $n\in\mathbb{Z}$ and each $k\in\mathbb{Z^{+}}$, there is an isomorphism \\
		$(A/\mathfrak{a}\otimes_{A}M)^{n}\cong (A/\mathfrak{a}^{k}\otimes_{A} M)^{n}$ . It follows that
		$\Lambda_{\mathfrak{a}}(M)^{n}:=(\underset{k}\varprojlim~ A/\mathfrak{a}^{k}\otimes_{A}M)^{n}\cong (A/\mathfrak{a}\otimes_{A}M)^{n}$ for all $n$.
		So, the complexes $\Lambda_{\mathfrak{a}}(M)$ and $A/\mathfrak{a}\otimes_{A} M$ are isomorphic.
		\item[$(3)\Rightarrow (4)$] Let $n\in\mathbb{Z}$. We have $A$-module isomorphisms
		$\mathfrak{a}\Lambda_{\mathfrak{a}}(M)^{n}\cong \mathfrak{a}(A/\mathfrak{a}\otimes_{A}M)^{n}\cong \mathfrak{a}(A/\mathfrak{a}\otimes_{A}M^{n})\cong 0$. So, $\mathfrak{a}\Lambda_{\mathfrak{a}}(M)\cong 0$. 
		\item[$(4)\Rightarrow (1)$] Since $\mathfrak{a}\Lambda_{\mathfrak{a}}(M)\cong 0$, the $A$-modules $ \mathfrak{a}\Lambda_{\mathfrak{a}}(M)^{n}\cong \mathfrak{a}\Lambda_{ \mathfrak{a}}(M^{n})=0$ for all $n$. It follows by \cite[Proposition 2.3]{David-Application-I} that $M^{n}$ is $\mathfrak{a}$-coreduced for each $n$. Therefore, $M\in \text{C}(A)_{\mathfrak{ a\text{-cor}}}$.
	\end{itemize}
\end{prf}
\begin{rem}\label{Rem: red=cor=any complex}
	Let $\mathfrak{ a}=\mathfrak{ a}^{2}$. We have {\normalfont $\text{C}(A)_{\mathfrak{ a\text{-red}}}=\text{C}(A)_{\mathfrak{ a\text{-cor}}}=\text{C}(A)$}.
\end{rem}

%\begin{cor}\label{Cor: product of cor is cor...}
%Let $\mathfrak{ a}$ be finitely generated. If {\normalfont $ \{X_{j}\}_{j\in\Lambda}\in \text{C}(A)_{\mathfrak{ a\text{-cor}}}$}, then so is $\underset {j\in\Lambda}\prod X_{j}$.
%\end{cor}
%\begin{prf}
%	First note that the functors $A/\mathfrak{ a}^{t}\otimes_{A}-, t\ge 1$ and $\Lambda_{ \mathfrak{a}}(-)$ commute with direct products \cite[2.2.7]{Peter-Schenzez}.  This together with Proposition \ref{sec2: Pro I-cor} yields
%	$\Lambda_{ \mathfrak{a}}(\underset{j\in\Lambda}{\prod}X_{j})=\underset{j\in\Lambda}{\prod}\Lambda_{ \mathfrak{a}}(X_{j}) \cong \underset{j\in\Lambda}\prod A/\mathfrak{ a}\otimes_{A}X_{j}\cong A/\mathfrak{ a}\otimes_{A}\underset{j\in\Lambda}\prod X_{j}.$
%\end{prf}
\begin{prop}\label{Prop: completion of I-cored}
	Let $\mathfrak{a}$ be a finitely generated ideal of $A$, and let {\normalfont $M\in\text{C}(A)$}. Then {\normalfont $M\in \text{C}(A)_{\mathfrak{ a\text{-cor}}}$} if and only if {\normalfont $\Lambda_{ \mathfrak{a}}(M)\in\text{C}(A)_{\mathfrak{ a\text{-cor}}}$}.
\end{prop}
\begin{prf}
	Let {\normalfont $M\in \text{C}(A)_{\mathfrak{ a\text{-cor}}}$}. By Proposition \ref{sec2: Pro I-cor}, $\mathfrak{ a}\Lambda_{ \mathfrak{a}}(M)=0$.  Since $\Lambda_\mathfrak{a}$ is an $A$-linear,
	$\mathfrak{a}\Lambda_{\mathfrak{a}}(\Lambda_{\mathfrak{a}}(M))\cong \Lambda_{ \mathfrak{a}}(\mathfrak{a}\Lambda_{\mathfrak{a}}(M))\cong 0.$
	So,  $\Lambda_{\mathfrak{a}}(M)\in \text{C}(A)_{\mathfrak{ a\text{-cor}}}$. Conversely, let $\Lambda_{\mathfrak{a}}(M)\in \text{C}(A)_{\mathfrak{ a\text{-cor}}}$. This together with the fact that $\Lambda_{ \mathfrak{a}}(-)$ is idempotent \cite[Section 3]{porta2014homology} and Proposition \ref{sec2: Pro I-cor},
	$\Lambda_{ \mathfrak{a}}(M)\cong \Lambda_{ \mathfrak{a}}(\Lambda_{ \mathfrak{a}}(M))\cong A/\mathfrak{ a}\otimes_{A}\Lambda_{ \mathfrak{a}}(M)\cong \Lambda_{ \mathfrak{a}}(M)/\mathfrak{ a}\Lambda_{ \mathfrak{a}}(M),$ which holds if and only if $\mathfrak{ a}\Lambda_{ \mathfrak{a}}(M)\cong 0$ if and only if $M\in\text{C}(A)_{\mathfrak{ a\text{-cor}}}$.
\end{prf}
\begin{cor}
	Let $\mathfrak{ a}$ be a finitely generated ideal of $A$, and {\normalfont $M,J\in\text{C}(A)$} with $J$ a complex of injective $A$-modules. 
	\begin{enumerate}
		\item [$(1)$]{\normalfont $\text{Hom}_{A}(M,J)\in \text{C}(A)_{\mathfrak{ a\text{-cor}}}$} if and only if {\normalfont $\text{Hom}_{A}(\Gamma_{ \mathfrak{a}}(M),J)\in \text{C}(A)_{\mathfrak{ a\text{-cor}}}$.} 
		\item[$(2)$]   $J$ is $\mathfrak{ a}$-coreduced if and only if {\normalfont $\text{Hom}_{A}(\Gamma_{ \mathfrak{a}}(A),J)$} is $\mathfrak{ a}$-coreduced.
	\end{enumerate}
\end{cor}
\begin{prf}
The proof of part $(1)$ is immediate by Proposition \ref{Prop: completion of I-cored}, since by \cite[Proposition 9.2.7]{Peter-Schenzez} {\normalfont $\text{Hom}_{A}(\Gamma_{ \mathfrak{a}}(M),J)\cong \Lambda_{ \mathfrak{a}}(\text{Hom}_{A}(M,J))$}. Take $M=A$ in part $(1)$ to prove part $(2)$.
\end{prf}
\begin{lem}\label{sec2: hom-red-cor module}
	Let $\mathfrak{a}$ be an ideal of $A$, $M$ an $\mathfrak{a}$-reduced module, and $N$ an $\mathfrak{a}$-coreduced module. Then for any module $X$ over $A$,
{\normalfont	\begin{itemize}
		\item[$(1)$]$\text{Hom}_{A}(X,M)\in \text{M}(A)_{\mathfrak{a}\text{-red}}$.
		\item[$(2)$]	$\text{Hom}_{A}(N,X)\in \text{M}(A)_{\mathfrak{a}\text{-cor}}$.
	\end{itemize}}
\end{lem}
\begin{prf}
	\begin{itemize}
		\item[$(1)$] Let  $\varphi\in \text{Hom}_{A}(X,M)$ and $a^{2}\varphi=0$ for all $a\in A$. It follows that $a^{2}\varphi(x)=0$ for all $a$ and for all $x\in X$. Since $\varphi(x)\in M$ and $M$ is $\mathfrak{a}$-reduced, 
		$a\varphi(x)=0$ for all $x$ and for all $a\in\mathfrak{ a}$. So, $a\varphi=0$ for any $a\in\mathfrak{a}$.  The conclusion is that the $A$-module $\text{Hom}_{A}(X,M)$ is  $\mathfrak{a}$-reduced by \cite[Definition 2.1]{kyomuhangi2020locally}.
		\item[$(2)$] Let $N$ be an $\mathfrak{a}$-coreduced $A$-module. In general, we have
		$\mathfrak{a}^{2}\text{Hom}_{A}(N,X)\subseteq \mathfrak{a}\text{Hom}_{A}(N,X)$ for any $A$-module $X$. We prove the reverse inclusion. Let $x\in \mathfrak{a}\text{Hom}_{A}(N,X)$. Then 
		$x=\sum_{i=1}^{k}a_{i}f_{i}(n)$ where $k\in\mathbb{Z^{+}}, a_{i}\in \mathfrak{a}, n\in N$ and $f_{i}\in \text{Hom}_{A}(N,X)$ for each $i=1,2,\dots, k$. Thus
		$	x=\sum_{i=1}^{k}f_{i}(a_{i}n)$ with $a_{i}n\in \mathfrak{a}N=\mathfrak{a}^{2}N$ since $N$ is $\mathfrak{a}$-coreduced. So, $a_{i}n\in \mathfrak{a}^{2}N$ for all $i$. Therefore, $a_{i}m=\sum_{j=1}^{t}b_{j}n_{j}$ where $b_{j}\in \mathfrak{a}^{2}$ and $n_{j}\in N$. It follows that
		$x=\sum_{i=1}^{k}f_{i}(\sum_{j=1}^{t}b_{j}n_{j})=\sum_{i=1}^{k} \sum_{j=1}^{t} b_{j}f_{i}(n_{j})= \sum_{j=1}^{t} b_{j}(\sum_{i=1}^{k}f_{i}(n_{j}))\in \mathfrak{a}^{2}\text{Hom}_{A}(N,X)$.
	\end{itemize}
\end{prf}
\begin{prop}\label{sec2: Prop Hom(-,M) and Hom(M,-)}
	Let $M\in {\normalfont\text{C}(A)_{\mathfrak{a}\text{-red}}}$ and {\normalfont$  N\in \text{C}(A)_{\mathfrak{a}\text{-cor}}$}. For any complex $X$,
	{\normalfont \begin{itemize}
		\item [$(1)$]$\text{ Hom}_{A}(X,M)\in {\normalfont\text{C}(A)_{\mathfrak{a}\text{-red}}} $,
		\item[ $(2)$] $\text{Hom}_{A}(N,X)\in {\normalfont  \text{C}(A)_{\mathfrak{a}\text{-cor}}}$.
	\end{itemize}}
\end{prop}	
\begin{prf}
	Let $X:=(X^{i},d^{i}_{X})_{i\in\mathbb{Z}}$ and $M:=(M^{j},d^{j}_{M})_{j\in\mathbb{Z}}$ be $A$-complexes with $M$ $\mathfrak{a}$-reduced. The $n^{th}$ position of the complex  $\text{Hom}_{A}(X,M)$ is the $A$-module {\normalfont $\underset{n=j-i}{\prod}\text{Hom}_{A}(X^{i},M^{j})$}. Since, by Lemma \ref{sec2: hom-red-cor module} $(1)$ $\text{Hom}_{A}(X^{i},M^{j})$ is $\mathfrak{a}$-reduced for each $j$ and the direct product of $\mathfrak{a}$-reduced modules is again $\mathfrak{a}$-reduced \cite[Proposition 2.4]{kyomuhangi2020locally}, it follows that the complex $\text{Hom}_{A}(X,M)$ is $\mathfrak{a}$-reduced. The proof of (2) is similar.
\end{prf}
%\begin{cor}
%	Let $\mathfrak{a}$ be an ideal of a ring $A$, and let $f: X\rightarrow Y$ be a quasi-isomorphism between $\text{K}$-injective (resp. $\text{K}$-projective) $A$-complexes.  	
%	The morphism
%{\normalfont
%$	\text{Hom}_{A}(f,Z): \text{Hom}_{A}(Y, Z)\rightarrow \text{Hom}_{A}(X, Z)$}
%	(resp.
%{\normalfont $	\text{Hom}_{A}(Z',f): \text{Hom}_{A}(Z', X)\rightarrow \text{Hom}_{A}(Z', Y)$}
%	is quasi-isomorphism of $\mathfrak{a}$-reduced (resp. $\mathfrak{a}$-coreduced) complexes for any $\mathfrak{a}$-reduced complex $Z$ and $\mathfrak{a}$-coreduced complex $Z'$.	
%	
%\end{cor}
%\begin{prf}
%	Let $f: X\rightarrow Y$ be a quasi-isomorphism between $K$-injective (resp. $K$-projective) $A$-complexes. Then by \cite[Corollary 4.4.11]{Peter-Schenzez} $f$ induces quasi-isomorphisms
%	$\text{Hom}_{A}(f,Z): \text{Hom}_{A}(Y, Z)\rightarrow \text{Hom}_{A}(X, Z)$
%	(resp. 
%$	\text{Hom}_{A}(Z',f): \text{Hom}_{A}(Z',X)\rightarrow \text{Hom}_{A}(Z', Y))$
%	for any $A$-complexes $Z$ and $Z'$. The conclusion is clear by Proposition \ref{sec2: Prop Hom(-,M) and Hom(M,-)}.
%\end{prf}
\begin{paragraph}\noi Proposition \ref{sec2: Prop Hom(-,M) and Hom(M,-)} says that the functors $\text{Hom}_{A}(M,-)$ (resp. $\text{Hom}_{A}(-,M)$) are endo functors on $\text{C}(A)_{\mathfrak{a}\text{-red}}$ (resp. $\text{C}(A)_{\mathfrak{a}\text{-cor}}$).
\end{paragraph}
\begin{lem}\label{Lem: cor tensor with cor is cor} Suppose that $M$ and $N$ are $A$-modules. If either $M$ or $N$ is $\mathfrak{ a}$-coreduced, then so is the $A$-module $M\otimes_{A}N$.
\end{lem}
\begin{prf}
	Suppose that either $aM= a^{2}M$ or $a^{2}N=aN$ for all $a\in\mathfrak{ a}$. Then either $aM\otimes_{A}N\cong a^{2}M\otimes_{A}N ~\text{or}~ M\otimes_{A}aN\cong M\otimes_{A}a^{2}N~\text{for all}~ a\in\mathfrak{ a}.$ In each case,  $a(M\otimes_{A}N)\cong a^{2}(M\otimes_{A}N)$ for every $a\in \mathfrak{ a}$. Therefore, $M\otimes_{A}N$ is $\mathfrak{ a}$-coreduced.
\end{prf}
\begin{prop}\label{prop: When tensor product of cored is cored} Let $\mathfrak{ a}$ be an ideal of $A$, and let {\normalfont $M,N\in\text{C}(A)$}. If either $M$ or $N$ is $\mathfrak{ a}$-coreduced, then so is $M\otimes_{A}N$. 
\end{prop}
\begin{prf}
	Since the element at the $n^{th}$-position of the complex $M\otimes_{A}N$ is defined by $(M\otimes_{A}N)^{n}=\underset{i+j=n}\bigoplus M^{i}\otimes_{A}N^{j}$ \cite[ 1.1.5]{Peter-Schenzez}, and $\mathfrak{ a}$-coreduced modules are closed under arbitrary direct sums, the conclusion follows by Lemma \ref{Lem: cor tensor with cor is cor}.
\end{prf}
\begin{prop}\label{Prop: leading to duality}
	Let $\mathfrak{ a}$ be an ideal of  $A$, and let {\normalfont $M,N\in \text{C}(A)$} where $N$ is an injective cogenerator of {\normalfont $\text{C}(A)$}.
	\begin{itemize}
		\item [$(1)$] The complex {\normalfont $M\in \text{C}(A)_{\mathfrak{ a\text{-cor}}}$} if and only if {\normalfont $\text{Hom}_{A}(M,N)\in \text{C}(A)_{\mathfrak{ a\text{-red}}}$}. 
		\item[$(2)$] Let $\mathfrak{ a}$ be finitely generated and {\normalfont $M\in \text{C}(A)_{\mathfrak{ a\text{-red}}}$}. The complex  {\normalfont $\text{Hom}_{A}(M,N)\in\text{C}(A)_{\mathfrak{ a\text{-cor}}}$}.
	\end{itemize} 
\end{prop}
\begin{prf}
	\begin{itemize}
		\item[$(1)$] The necessity part is proved in Proposition \ref{sec2: Prop Hom(-,M) and Hom(M,-)} $(2)$. We prove only  the sufficiency. Suppose that {\normalfont $\text{Hom}_{A}(M,N)\in\text{C}(A)_{\mathfrak{ a\text{-red}}}$}. This together with the Hom-Tensor adjunction gives
		{\normalfont $\text{Hom}_{A}(A/\mathfrak{ a}\otimes_{A}M,N)\cong  \text{Hom}_{A}(A/\mathfrak{ a},\text{Hom}_{A}(M,N))\cong 
			\text{Hom}_{A}(A/\mathfrak{ a}^{2},\text{Hom}_{A}(M,N))\cong \text{Hom}_{A}(A/\mathfrak{ a}^{2}\otimes_{A}M,N)$}.	
		But since $N$ is an injective cogenerator, the functor $\text{Hom}_{A}(-,N)$ reflects isomorphisms. So, $ A/\mathfrak{ a}\otimes_{A}M \cong A/\mathfrak{ a}^{2}\otimes_{A}M$. It follows, by Proposition \ref{sec2: Pro I-cor}, that $M\in\text{C}(A)_{\mathfrak{ a\text{-cor}}}$.
		\item[$(2)$]By \cite[Lemma 1.4.6]{Peter-Schenzez} and by Proposition \ref{sec2: Pro-I-red}, the isomorphisms $A/\mathfrak{ a}\otimes_{A}\text{Hom}_{A}(M,N)\cong \text{Hom}_{A}(\text{Hom}_{A}(A/\mathfrak{ a},M), N)\cong
		\text{Hom}_{A}(\text{Hom}_{A}(A/\mathfrak{ a}^{2},M), N)\cong A/\mathfrak{ a}^{2}\otimes_{A} \text{Hom}_{A}(M,N)$ are evident. The conclusion $\text{Hom}_{A}(M,N)\in\text{C}(A)_{\mathfrak{ a\text{-cor}}}$ holds by Proposition \ref{sec2: Pro I-cor}.
	\end{itemize}
\end{prf}
\begin{paragraph}\noi
	Denote by $(.)^{\vee}:= \text{Hom}_{A}(-,E)$ the general Matlis dual endo-functor on {\normalfont $\text{C}(A)$}, where $E$ is an injective cogenerator of {\normalfont $\text{C}(A)$}. 
\end{paragraph}
\begin{cor}\label{Cor: X is cor iff its dual is red}
	Let $\mathfrak{ a}$ be an ideal of  $A$, and let {\normalfont $M\in\text{C}(A)$}. Then {\normalfont $M\in\text{C}(A)_{\mathfrak{ a\text{-cor}}}$}  if and only if {\normalfont $M^{\vee}\in\text{C}(A)_{\mathfrak{ a\text{-red}}}$}. If {\normalfont $M\in\text{C}(A)_{\mathfrak{ a\text{-red}}}$}, then {\normalfont $M^{\vee}\in\text{C}(A)_{\mathfrak{ a\text{-cor}}}$} provided $\mathfrak{ a}$ is finitely generated.
\end{cor}
\begin{prf}	By taking $N:= E$, Proposition \ref{Prop: leading to duality} $(1)$ (resp. $(2)$) provides  proof of the first (resp. the second) statement.	
\end{prf} 
\begin{cor}
	Let $\mathfrak{ a}$ be an ideal of a Noetherian ring $A$, and let {\normalfont $M,N\in \text{C}(A)$} with {\normalfont $M \in \text{C}(A)_{\mathfrak{ a\text{-red}}}$}. 
	\begin{enumerate}
		\item[$(1)$] {\normalfont$\Lambda_{ \mathfrak{a}}(M^{\vee})\in  \text{C}(A)_{\mathfrak{ a\text{-cor}}}$}.
		\item[$(2)$] {\normalfont $\text{Hom}_{A}(\Lambda_{ \mathfrak{a}}(M),\Lambda_{ \mathfrak{a}}(N))\in\text{C}(A)_{\mathfrak{ a\text{-cor}}}$}.
		\item[$(3)$] If, as a complex, $A$ is an $\mathfrak{ a}$-reduced ring which $\mathfrak{ a}$-adically complete, then {\normalfont $\Lambda_{ \mathfrak{a}}(N)\in\text{C}(A)_{\mathfrak{ a\text{-cor}}}$}.
	\end{enumerate}  
\end{cor}
\begin{prf}
	Since {\normalfont $M\in  \text{C}(A)_{\mathfrak{ a\text{-red}}}$}, $M^{\vee}\in \text{C}(A)_{\mathfrak{ a\text{-cor}}}$ [Proposition \ref{Prop: leading to duality} $(2)$], and the conclusion $\Lambda_{ \mathfrak{a}}(M^{\vee})\in \text{C}(A)_{\mathfrak{ a\text{-cor}}}$ holds by Proposition \ref{Prop: completion of I-cored}. The proof of part $(2)$ is immediate by Proposition \ref{Prop: leading to duality} $(2)$, since {\normalfont $\text{Hom}_{A}(\Lambda_{ \mathfrak{a}}(M),\Lambda_{ \mathfrak{a}}(N))\cong \text{Hom}_{A}( M,\Lambda_{ \mathfrak{a}}(N))$} for any complex $N$ \cite[Lemma 7.2.2]{Peter-Schenzez}. By taking $M=A$ in part $(2)$, proof of part $(3)$ follows.
\end{prf}
\begin{paragraph}\noi Proposition \ref{Prop: Properties} summarises some closure properties. The proof can be adapted from \cite{kyomuhangi2020locally,David-Application-I} studied in the case of {\normalfont $\text{M}(A)$}.
	\end{paragraph}
\begin{prop} \label{Prop: Properties}
	Let $\mathfrak{a}$ be an ideal of  $A$.  Let $\{M_{j}\}_{j\in\Lambda}$ (resp. $\{N_{k}\}_{k\in\Lambda}$) be a family of $\mathfrak{ a}$-reduced complexes (resp. $\mathfrak{ a}$-coreduced complexes) where $\Lambda$ is an index set. The following statements hold. 
	\begin{itemize}
		\item [$(1)$] {\normalfont $\prod_{j\in\Lambda}^{} M_{j} \in\text{C}(A)_{\mathfrak{ a\text{-red}}}$} and  {\normalfont $\bigoplus_{k\in\Lambda} N_{k} \in\text{C}(A)_{\mathfrak{ a\text{-cor}}}$}. 
		%\item [$(2)$] $\bigoplus_{k\in\Lambda}^{} N_{k} $ is  $\mathfrak{ a}$-coreduced with respect to $M$.
		\item [$(2)$] If $\{M_{j}\}_{j\in\Lambda}$ is an inverse system of $A$-complexes, then {\normalfont $\underset{j}\varprojlim M_{j}\in\text{C}(A)_{\mathfrak{ a\text{-red}}}$}. 
		
		\item[$(3)$] If $\{N_{k}\}_{k\in\Lambda}$ is a direct system of $A$-complexes, then {\normalfont $\underset{k}\varinjlim N_{k}\in\text{C}(A)_{\mathfrak{ a\text{-cor}}}$}. 
		\item [$(4)$] {\normalfont $\text{C}(A)_{\mathfrak{ a\text{-red}}}$} (resp. {\normalfont $\text{C}(A)_{\mathfrak{ a\text{-cor}}}$}) is closed under submodules (resp. under quotients).
		\end{itemize} 
		\end{prop}
%	\begin{paragraph}\noi Proposition \ref{Prop: Properties}, in particular, says that
%		\begin{enumerate}
%		\item[$(a)$] \textit{the right soft truncation} of an $\mathfrak{ a}$-reduced cochain complex $M$ at spot $t$ which is given by  $\tau_{t]}~:~\dots\rightarrow M^{t-2}\rightarrow M^{t-1}\rightarrow \text{Z}^{t}(M)\rightarrow 0$, where $\text{Z}^{t}(M):=\text{Ker}(d^{t}_{M})$ is $\mathfrak{ a}$-reduced, and   \textit{the left soft truncation} of an $\mathfrak{ a}$-coreduced chain complex $M$ at spot $t$ given by $\tau_{[t}: 0\rightarrow M_{t}/\text{B}_{t}(M)\rightarrow M_{t-1}\rightarrow M_{t-2}\rightarrow\dots$, where $\text{B}_{t}(M):= \text{Im}(d^{M}_{t+1})$, is $\mathfrak{ a}$-coreduced.
%		\item[$(b)$] If $M\in\text{C}(A)_{\mathfrak{ a\text{-cor}}}$, then so is its ($i^{th}$-) homology $\text{H}_{i}(M)$.
%		 \end{enumerate}
%	\end{paragraph}

\begin{paragraph}\noi Recall that a category $\mathcal{A}$ is \textit{complete} if $\underset{}\varprojlim M_{j}$ exists in $\mathcal{A}$ for every inverse system $\{M_{j}\}_{j\in\Lambda}$ in $\mathcal{A}$. Dually,  $\mathcal{A}$ is \textit{cocomplete} if $\underset{}\varinjlim N_{j}$ exists in $\mathcal{A}$ for every direct system $\{N_{j}\}_{j\in\Lambda}$ in $\mathcal{A}$; see \cite[Definition 5.27]{J.Rotman- Introduction hom.algebra}. Moreover, a class $\mathscr{C}$ of objects is called a \textit{pretorsion class} if it is closed under quotient objects and coproducts, and is a \textit{pretorsion-free class} if it is closed under subobjects and products \cite[page 137]{Stenstrom--rings of quotients}. By Proposition \ref{Prop: Properties} it follows that
\end{paragraph}
\begin{prop}\label{prop: complete vs cocomplete}
	The full subcategory  {\normalfont $\text{C}(A)_{\mathfrak{ a}\text{-red}}$} (resp. {\normalfont $\text{C}(A)_{\mathfrak{ a}\text{-cor}}$}) of {\normalfont $\text{C}(A)$} is complete and forms a pretorsion-free class (resp. cocomplete and forms a pretorsion class). 
\end{prop}

%\begin{paragraph}\noi Suppose that $\text{C}(A)_{\mathfrak{ a\text{-red}}}$ (resp. $\text{C}(A)_{\mathfrak{ a\text{-cor}}}$)  have enough injectives (resp. enough projectives). If $M\rightarrow J$ is an injective resolution of $M$ with $M\in\text{C}(A)_{\mathfrak{ a\text{-red}}}$, then $J\in\text{C}(A)_{\mathfrak{ a\text{-red}}}$.  If $M\rightarrow J$ is an injective resolution of $M$ with $M\in\text{C}(A)_{\mathfrak{ a\text{-red}}}$, then $J\in\text{C}(A)_{\mathfrak{ a\text{-red}}}$. Dually, if $P\rightarrow M$ is a projective resolution of $M$ such that $M\in\text{C}(A)_{\mathfrak{ a\text{-cor}}}$, then $P\in\text{C}(A)_{\mathfrak{ a\text{-cor}}}$. \newline
%	Give an \textbf{example} to show that injetive resolutions (resp. projective resolutions) of $M\in\text{C}(A)_{\mathfrak{ a\text{-red}}}$ (resp. $M\in\text{C}(A)_{\mathfrak{ a\text{-cor}}}$) need not be in $\text{C}(A)_{\mathfrak{ a\text{-red}}}$ (resp. in $\text{C}(A)_{\mathfrak{ a\text{-cor}}}$).
%\end{paragraph}

%\section{Associated exact sequences}
%\begin{paragraph}\noi The full subcategory $\text{C}(A)_{\mathfrak{a}\text{-red}}$ (resp. $\text{C}(A)_{\mathfrak{a}\text{-cor}}$  ) of $\text{C}(A)$ is not closed under taking quotients (resp. sub-objects). The short exact sequences in the corresponding full subcategories are not preserved under a given, even, exact functor in general. In this section, we explore some exact functors which preserve short exact sequences of
%	$\mathfrak{ a}$-reduced complexes and $\mathfrak{ a} $-coreduced complexes.
%	\end{paragraph}

\begin{paragraph}\noi 
Let $f:M\rightarrow N$ be a morphism in $\text{C}(A)$. The \textit{cone} of $f$ is a chain complex $C_{f}:=(C_{f}^{n}, d_{C_{f}}^{n})$, where $C_{f}^{n}:= M^{n+1}\bigoplus N^{n} $, and the differentials $d_{C_{f}}^{n}: C_{f}^{n}\rightarrow C_{f}^{n+1}$ are given by 

\[
d_{C_{f}}^{n}=\begin{bmatrix}
	-d_{M}^{n+1}       & 0 \\
	f^{n+1}       & d_{N}^{n}  \\
\end{bmatrix}.
\]
The \textit{fibre} $F_{f}$ of $f$ is defined by $F_{f}:= C^{\left[-1\right]}_{f}$, i.e., $F_{f}:=(F_{f}^{n}, d_{F_{f}}^{n})$ where $F_{f}^{n}:= M^{n}\bigoplus N^{n-1}$ and 

\[
d_{F_{f}}^{n}=\begin{bmatrix}
	d_{M}^{n}       & 0 \\
	f^{n}       & -d_{N}^{n-1}  \\
\end{bmatrix}.
\]
\end{paragraph}
\begin{paragraph}\noi
Associated to any direct system $\mathscr{D}:=\{\sigma_{n,n+1}:M^{n}\rightarrow M^{n+1}| n\in\Z^{+}\}$ over $\Z^{+}$ of complexes $M^{n}$ is a morphism$$ \phi_\mathscr{D} : \bigoplus M^{n}\rightarrow \bigoplus M^{n},~ x\mapsto x-\sigma_{n,n+1}(x)$$ for $x\in M^{n}$ \cite[1.3.2]{Peter-Schenzez}. If $\mathscr{B}:= \{\rho_{n,n+1}:N^{n+1}\rightarrow N^{n}| n\in\Z^{+}\}$ is an inverse system of complexes, then there is a morphism associated to $\mathscr{B}$
$$\psi_\mathscr{B} : \prod N^{n}\rightarrow \prod N^{n},~ (y_{n})\mapsto (y_{n}-\rho_{n,n+1}(y_{n+1}))$$ for $(y_{n})\in\prod N^{n}$ \cite[1.2.2]{Peter-Schenzez}. Following \cite{greenlees1992derived, Peter-Schenzez}, the \textit{microscope} of $\mathscr{B}$ is defined as the fiber of $\psi_{\mathscr{B}}$, i.e., $\text{Mic}(\mathscr{B}):= F_{\psi_{\mathscr{B}}}$, and the \textit{telescope} of $\mathscr{D}$ is the cone of $\phi_\mathscr{D}$. So, $\text{Tel}(\mathscr{D}):= C_{\phi_\mathscr{D}}$.
	\end{paragraph}

\begin{prop}\label{Prop: Telescope-cored }
	Let $\mathscr{D}_{i},i=1,2,3$ be direct systems of complexes over $\Z^{+}$. 
	\begin{enumerate}
		\item [$(1)$] If {\normalfont $\mathscr{D}_{i}\in\text{C}(A)_{\mathfrak{ a\text{-cor}}}$} for each $i$, then so are {\normalfont $\text{Tel}(\mathscr{D}_{i})$} for each $i$.
		
		\item[$(2)$] Let $0\rightarrow \mathscr{D}_{1}\rightarrow \mathscr{D}_{2}\rightarrow \mathscr{D}_{3}\rightarrow 0$ be a short exact sequence of direct sytems of $\mathfrak{ a}$-coreduced complexes over $\Z^{+}$. Then the induced sequence {\normalfont $$0\rightarrow \text{Tel}(\mathscr{D}_{1})\rightarrow \text{Tel}(\mathscr{D}_{2})\rightarrow \text{Tel}(\mathscr{D}_{3})\rightarrow 0$$} is a short exact sequence of $\mathfrak{ a}$-coreduced complexes.
	\end{enumerate}
\end{prop}
\begin{prf}
	Let $(\mathscr{D}_{i}, \phi_{\mathscr{D}_{i}}), i=1,2,3$ be a pair consisting of direct sytems of $\mathfrak{ a}$-coreduced complexes and the associated morphisms. Since $\mathfrak{ a}$-coreduced complexes are closed under taking direct sums [Proposition \ref{Prop: Properties}], the cone $C_{\phi_{\mathscr{D}_{i}}}\in \text{C}(A)_{\mathfrak{ a\text{-cor}}}$ for each $i$. So, $\text{Tel}(\mathscr{D}_{i})\in \text{C}(A)_{\mathfrak{ a\text{-cor}}}$ for each $i$.  This proves part $(1)$. The proof of part $(2)$ follows by part $(1)$ and \cite[Lemma 4.3.4]{Peter-Schenzez}.
\end{prf}
\begin{prop}\label{Prop: Microscope }
	Let $\mathscr{B}_{i},i=1,2,3$ be inverse systems of complexes over $\Z^{+}$. 
	\begin{enumerate}
		\item [$(1)$] If {\normalfont $\mathscr{B}_{i}\in\text{C}(A)_{\mathfrak{ a\text{-red}}}$} for each $i$, then so are {\normalfont $\text{Mic}(\mathscr{B}_{i})$} for each $i$.
		
		\item[$(2)$] Let $0\rightarrow \mathscr{B}_{1}\rightarrow \mathscr{B}_{2}\rightarrow \mathscr{B}_{3}\rightarrow 0$ be a short exact sequence of inverse sytems of $\mathfrak{ a}$-reduced complexes over $\Z^{+}$. Then the induced sequence {\normalfont $$0\rightarrow \text{Mic}(\mathscr{B}_{1})\rightarrow \text{Mic}(\mathscr{B}_{2})\rightarrow \text{Mic}(\mathscr{B}_{3})\rightarrow 0$$} is a short exact sequence of $\mathfrak{ a}$-reduced complexes.
	\end{enumerate}
\end{prop}
\begin{prf}
	Let $(\mathscr{B}_{i}, \psi_{\mathscr{B}_{i}}), i=1,2,3$ be a pair consisting of inverse systems of $\mathfrak{ a}$-reduced complexes and the associated morphisms. Since $\mathfrak{ a}$-reduced complexes are closed under direct sums  by virtue of \cite[Proposition 3.4]{kyomuhangi2020locally}, the fiber $F_{\psi_{\mathscr{B}_{i}}}$, which is the cone $C_{\psi_{\mathscr{B}_{i}}}^{\left[-1\right]}$, belongs to $\text{C}(A)_{\mathfrak{ a\text{-red}}}$ for each $i = 1, 2, 3$. It follows that $\text{Mic}(\mathscr{B_{i}}) \in \text{C}(A)_{\mathfrak{ a\text{-red}}}$. This proves part $(1)$. The proof of part $(2)$ holds by \cite[Lemma 4.2.5]{Peter-Schenzez}
	and part (1).
\end{prf}

\section{Results in $\text{C}(A)$} \label{ Section: C(A) staff}
\begin{paragraph}\noi The purpose of this section is to prove versions of the GM duality and the MGM equivalence in $\text{C}(A)$. We denote by $\mathfrak{C}(A)$ the full subcategory $ \text{C}(A)_{\mathfrak{a}\text{-red}} ~\cap~ \text{C}(A)_{\mathfrak{a}\text{-cor}}$ of $\text{C}(A)$.
	\end{paragraph}
\begin{lem}\label{sec3: Idempotent-composition}
	Let $\mathfrak{a}$ be an ideal of a ring $A$, and let $M$ be a complex of $A$-modules.  
	\begin{itemize}
		\item [$(1)$] The functors {\normalfont $$A/\mathfrak{a}~\otimes_{A}-,\text{ Hom} _{A}(A/\mathfrak{a},-): \text{C}(A)\rightarrow \mathfrak{C}(A)$$} are idempotent.
		\item[$(2)$] For any complex $M$,	
		{\normalfont $
		A/\mathfrak{a} ~\otimes_{A}\text{	Hom}_{A}(A/\mathfrak{a},M)\cong \text{Hom}_{A}(A/\mathfrak{a},M)$ }\text{and}
		{\normalfont $$\text{Hom}_{A}(A/\mathfrak{a}, A/\mathfrak{a}\otimes_{A} M)\cong A/\mathfrak{a}\otimes_{A} M.$$}
%		\item[$(3)$] For any $A$-complex $M$, {\normalfont $\text{Hom}_{A}(A/\mathfrak{a},M)\in \text{C}(A)_{\mathfrak{a\text{-tor}}}$ and $A/\mathfrak{a}~\otimes_{A} M\in \text{C}(A)_{\mathfrak{a\text{-com}}}$}.
	\end{itemize}
\end{lem}
\begin{prf}
\begin{enumerate}
	\item[$(1)$] Let $ M=(M^{n})_{n\in\mathbb{Z}}\in \text{C}(A)$. Set $L=A/\mathfrak{a}\otimes_{A}M$. Since $A/\mathfrak{a}$ is concentrated in degree $0$, by Lemma \ref{sec2: Hom-tensor Lem}
	\begin{equation*}
	(A/\mathfrak{a}~ \otimes_{A} L)^{n}=A/\mathfrak{a}~\otimes_{A} L^{n}= A/\mathfrak{a}~\otimes_{A} (A/\mathfrak{a}~\otimes_{A}~ M)^{n}\cong A/\mathfrak{a}~\otimes_{A} (A/\mathfrak{a}~\otimes_{A}~ M^{n})
	\end{equation*}
	\begin{equation*}
	\cong A/\mathfrak{a}~\otimes_{A} M^{n} \cong (A/\mathfrak{a}~\otimes_{A}~ M)^{n}. 
	\end{equation*}
	
	Thus
	\begin{equation*}
	(A/\mathfrak{a}\otimes_{A}(A/\mathfrak{a}\otimes_{A}M))^{n}=(A/\mathfrak{a}\otimes_{A}M)^{n}
	\end{equation*}
	for all $n$.
	It follows that $A/\mathfrak{a}\otimes_{A}(A/\mathfrak{a}\otimes_{A}M)\cong A/\mathfrak{a}\otimes_{A}M$. Therefore, the functor $A/\mathfrak{a}\otimes_{A}-$ is idempotent on $\text{C}(A)$.
	For the second idempotence it suffices to show that for any $n\in\mathbb{Z}$,
	\begin{equation*}
	\text{Hom}_{A}(A/\mathfrak{a},\text{Hom}_{A}(A/\mathfrak{a},M)^{n}= \text{Hom}_{A}(A/\mathfrak{a},M)^{n}.
	\end{equation*}
	Set $N:= \text{Hom}_{A}(A/\mathfrak{a},M)$. By Lemma \ref{sec2: Hom-tensor Lem}, $ N^{n}:= \text{Hom}_{A}(A/\mathfrak{a},M^{n}).$
	
	So,
	$\text{Hom}_{A}(A/\mathfrak{a},\text{Hom}_{A}(A/\mathfrak{a},M))^{n}= \text{Hom}_{A}(A/\mathfrak{a},N)^{n}=\text{Hom}_{A}(A/\mathfrak{a},N^{n}) 
	=\text{Hom}_{A}(A/\mathfrak{a},\text{Hom}_{A}(A/\mathfrak{a},M^{n})).$
	By \cite[Proposition 2.9]{David-Application-I}
	\begin{equation*}
	\text{Hom}_{A}(A/\mathfrak{a},\text{Hom}_{A}(A/\mathfrak{a},M^{n}))\cong \text{Hom}_{A}(A/\mathfrak{a},M^{n})= \text{Hom}_{A}(A/\mathfrak{a},M)^{n}.
	\end{equation*}
	Therefore for any $A$-complex $M$,
	\begin{equation*}
	\text{Hom}_{A}(A/\mathfrak{a},\text{Hom}_{A}(A/\mathfrak{a},M))\cong \text{Hom}_{A}(A/\mathfrak{a},M).
	\end{equation*}
	For the remaining part, let $M$ be an $A$-complex and $n\in\mathbb{Z}$.
	By \cite[Proposition 2.9]{David-Application-I} the $A$-modules $ \text{Hom}_{A}(A/\mathfrak{a}, M^{n})$ and $A/\mathfrak{a}\otimes_{A} M^{n}$ are both $\mathfrak{ a}$-reduced and $\mathfrak{ a}$-coreduced for each $n$. By Lemma \ref{sec2: Hom-tensor Lem} it follows that $\text{Hom}_{A}(A/\mathfrak{a}, M)^{n}$ and 
	$(A/\mathfrak{a}\otimes_{A} M)^{n}$ are $\mathfrak{ a}$-reduced and $\mathfrak{ a}$-coreduced for all $n$. This shows that the $A$-complexes
	$\text{Hom}_{A}(A/\mathfrak{a}, M)$ and $A/\mathfrak{a}~\otimes_{A} M $ are both $\mathfrak{ a}$-reduced and $\mathfrak{ a}$-coreduced. So, they are contained in $\mathfrak{C}(A)$.	 
	\item[$(2)$]
	\begin{itemize}
		\item[i.]$	
		(A/\mathfrak{a} ~\otimes_{A} \text{Hom}_{A}(A/\mathfrak{a},M))^{n}= A/\mathfrak{a} ~\otimes_{A} (\text{Hom}_{A}(A/\mathfrak{a},M))^{n}
		= A/\mathfrak{a} ~\otimes_{A} \text{Hom}_{A}(A/\mathfrak{a},M^{n})$.	However, by \cite[Proposition 2.9]{David-Application-I}
		$$A/\mathfrak{a} ~\otimes_{A} (\text{Hom}_{A}(A/\mathfrak{a},M^{n})\cong \text{Hom}_{A}(A/\mathfrak{a},M^{n})~ \text{for all} ~n. $$ So,
		$$(A/\mathfrak{a} ~\otimes_{A} \text{Hom}_{A}(A/\mathfrak{a},M))^{n}\cong \text{Hom}_{A}(A/\mathfrak{a},M^{n})\cong \text{Hom}_{A}(A/\mathfrak{a},M)^{n}.$$
		Therefore for any complex $M$,	
		$ A/\mathfrak{a} ~\otimes_{A}	\text{Hom}_{A}(A/\mathfrak{a},M)\cong \text{Hom}_{A}(A/\mathfrak{a},M).$
		\item[ii.] We show that $$\text{Hom}_{A}(A/\mathfrak{a}, A/\mathfrak{a}\otimes_{A}M)^{n}=(A/\mathfrak{a}\otimes_{A}M)^{n}~ \text{for all}~ n.$$	$\text{Hom}_{A}(A/\mathfrak{a}, A/\mathfrak{a}~\otimes_{A}M)^{n}= \text{Hom}_{A}(A/\mathfrak{a},( A/\mathfrak{a}~\otimes_{A}M)^{{n}})
		=$ $$\text{Hom}_{A}(A/\mathfrak{a}, A/\mathfrak{a}~\otimes_{A}M^{{n}})\cong A/\mathfrak{a}~\otimes_{A}M^{{n}}~\text{by \cite[Proposition 2.9]{David-Application-I} for all $n$}.$$ So, there is an isomorphism of $A$-complexes $\text{Hom}_{A}(A/\mathfrak{a}, A/\mathfrak{a}\otimes_{A}M)\cong A/\mathfrak{a}\otimes_{A}M.$	
	\end{itemize}
%	\item[$(3)$] For an $A$-complex $M$ and $n\in\mathbb{Z}$,
%	$\text{Hom}_{A}(A/\mathfrak{a},M)^{n}= \text{Hom}_{A}(A/\mathfrak{a},M^{n})$  and $(A/\mathfrak{a}\otimes_{A}M)^{n}= A/\mathfrak{a}~\otimes_{A}M^{n}$ [Lemma \ref{sec2: Hom-tensor Lem}].
%	By \cite[Proposition 2.6]{David-Application-I}  the $A$-modules $\text{Hom}_{A}(A/\mathfrak{a},M^{n})$ (resp.$~A/\mathfrak{a}~\otimes_{A}M^{n}$) are $\mathfrak{a}$-torsion (resp. $\mathfrak{a}$-complete $A$-modules) for all $n$. Therefore the $A$-complexes $\text{Hom}_{A}(A/\mathfrak{a},M)$
%	(resp. $A/\mathfrak{a}~\otimes_{A}M$) are $\mathfrak{a}$-torsion (resp. $\mathfrak{a}$-adically complete)
%	complexes.
\end{enumerate}
\end{prf}
	
\begin{thm}[The version of GM duality in $\text{C}(A)$]\label{sec3: GM Duality-Thm}
	Let $A$ be a ring and $\mathfrak{a}$ an ideal in it.
	\begin{itemize}
		\item [$(1)$] The functor
		{\normalfont\begin{equation*}
			\Gamma_{\mathfrak{a}}: \text{C}(A)_{\mathfrak{a}\text{-red}}\rightarrow \text{C}(A)_{\mathfrak{a}\text{-cor}}, ~ M\mapsto \Gamma_{\mathfrak{a}}(M)
			\end{equation*}}
		is idempotent and {\normalfont $\Gamma_{\mathfrak{a}}(M)\cong \text{Hom}_{A}(A/\mathfrak{a},M)$}.
		\item[$(2)$] The functor 
		{\normalfont \begin{equation*}
			\Lambda_{\mathfrak{a}}: \text{C}(A)_{\mathfrak{a}\text{-cor}}\rightarrow \text{C}(A)_{\mathfrak{a}\text{-red}}, ~ M\mapsto \Lambda_{\mathfrak{a}}(M)
			\end{equation*}}
		is idempotent and {\normalfont $\Lambda_{\mathfrak{a}}(M)\cong A/\mathfrak{a}\otimes_{A} M$}.
		\item[$(3)$] For any $ M\in {\normalfont \text{C}(A)_{\mathfrak{a}\text{-cor}}} $ and any {\normalfont $ N\in \text{C}(A)_{\mathfrak{a}\text{-red}}$}, 
		{\normalfont	\begin{equation*}
			\text{Hom}_{A}( \Lambda_{\mathfrak{a}}(M),N)\cong 	\text{Hom}_{A}( M,\Gamma_{\mathfrak{a}} (N)).
			\end{equation*}}
	\end{itemize}
\end{thm}
\begin{prf}
	\begin{itemize}
		\item [$(1)$] For any $\mathfrak{a}$-reduced module $N$, $\Gamma_{\mathfrak{a}}(N)= \underset{k}\varinjlim~ {\text{Hom}_{A}(A/\mathfrak{a}^{k},N)}\cong \text{Hom}_{A}(A/\mathfrak{a},N)$ 
		is an $\mathfrak{ a}$-coreduced module \cite[Proposition 2.9]{David-Application-I}. Let $M$ be  an $\mathfrak{ a}$-reduced complex. Applying the functor $\Gamma_{\mathfrak{a}}$ levelwise to $M$ one gets $\Gamma_{\mathfrak{a}}(M)\in \text{C}(A)_{\mathfrak{ a\text{-cor}}}$. By Proposition \ref{sec2: Pro-I-red}, $\Gamma_{\mathfrak{a}}(M)\cong \text{Hom}_{A}(A/\mathfrak{a},M)$. The functor
		$ \text{Hom}_{A}(A/\mathfrak{a},-)$ is idempotent by Lemma \ref{sec3: Idempotent-composition}.
		\item[$(2)$] 
		For any $\mathfrak{a}$-coreduced complex $M$, $\Lambda_{\mathfrak{a}}(M)\cong A/\mathfrak{a}~\otimes_{A}M$ by Proposition \ref{sec2: Pro I-cor}. Since for any $\mathfrak{a}$-coreduced module $N$,  $\Lambda_{\mathfrak{a}}(N) = \underset{k}\varprojlim~ {N/\mathfrak{a}^{k}N}$ is an $\mathfrak{ a}$-reduced $A$- module \cite[Proposition 2.9]{David-Application-I},  levelwise application extends $\Lambda_{ \mathfrak{a}}$ to a complex of $\mathfrak{ a}$-reduced modules. Thus if $M$ is an $\mathfrak{ a}$-coreduced complex, then $\Lambda_{\mathfrak{a}}(M)\in \text{C}(A)_{\mathfrak{ a\text{-red}}}$. The idempotence of $\Lambda_{ \mathfrak{a}}(-)\cong A/\mathfrak{a}\otimes_{A} -$ follows by Lemma \ref{sec3: Idempotent-composition}.
		\item[$(3)$] The functor $ A/\mathfrak{a}\otimes_{A}-$ is left adjoint to $\text{Hom}_{A}(A/\mathfrak{a},-)$ in the category of $A$-modules. Levelwise application extends the functors to the category $\text{C}(A)$. Also for any $M\in\text{C}(A)_{\mathfrak{a}\text{-red}}$, $\Gamma_{\mathfrak{a}}(M)\cong \text{Hom}_{A}(A/\mathfrak{a},M)$ and for any $ N\in \text{C}(A)_{\mathfrak{a}\text{-cor}}$, $  A/\mathfrak{a}~\otimes_{A}N\cong \Lambda_{\mathfrak{a}}(N) $. By uniqueness of adjoints, $\Lambda_{ \mathfrak{a}}$ is left adjoint to $\Gamma_{ \mathfrak{a}}$. 
	\end{itemize}
\end{prf}
\begin{cor}\label{Lem: Exactness of the functors}
	{\normalfont Let $\mathfrak{ a}$ be an ideal of $A$.
		\begin{itemize}
			\item [$(1)$] The functor $\Gamma_{ \mathfrak{a}}(-): \text{C}(A)_{\mathfrak{a}\text{-red}}\rightarrow \text{C}(A)_{\mathfrak{a}\text{-cor}}$ is left exact. 
			\item [$(2)$] The functor $\Lambda_{ \mathfrak{a}}(-): \text{C}(A)_{\mathfrak{a}\text{-cor}}\rightarrow \text{C}(A)_{\mathfrak{a}\text{-red}}$ is right exact. 
	\end{itemize}}
\end{cor}
\begin{prf}
	This immediately follows from Theorem \ref{sec3: GM Duality-Thm} and by \cite[Theorem 2.6.1]{weibel1995introduction}.
\end{prf}
\begin{paragraph}\noi
	Denote by $\textbf{Set}$ the category of all sets. Recall that a functor $G: \mathscr{C}\rightarrow \textbf{Set}$ is \textit{representable} if $G$ is naturally isomorphic to $\text{Hom}_{R}(M,-)$ for some $M\in \mathscr{C}$.
\end{paragraph}
\begin{prop}\label{prop: representablity of gen. Gamma functor}
	Let $\mathfrak{ a}$ be an ideal of  $A$. A representable functor {\normalfont $\Gamma_{ \mathfrak{a}}(-): \text{C}(A)_{\mathfrak{a}\text{-red}}\rightarrow \text{C}(A)_{\mathfrak{a}\text{-cor}}$}  preserves all limits. In particular, {\normalfont $\Gamma_{ \mathfrak{a}}(-)$} preserves inverse limits, coproducts, pullbacks, equalizers, and  terminal objects. 
\end{prop}
\begin{prf}
	Representablity of $\Gamma_{ \mathfrak{a}}(-)$ follows, since {\normalfont $A/\mathfrak{ a}\in \text{C}(A)_{\mathfrak{ a\text{-red}}}$}  and $\Gamma_{ \mathfrak{a}}(M)\cong \text{Hom}_{A}(A/\mathfrak{ a},M)$ for any {\normalfont $M\in\text{C}(A)_{\mathfrak{ a\text{-red}}}$}. Moreover, since $\Gamma_{ \mathfrak{a}}(-)$ is right adjoint to $\Lambda_{ \mathfrak{a}}(-)$ [Theorem \ref{sec3: GM Duality-Thm}], it preserves limits by \cite[Theorem 2.6.10]{weibel1995introduction}. 
\end{prf}
\begin{paragraph}\noi  We get Proposition \ref{Prop: Lambda_I -left adjoint} by a similar proof. 
	\end{paragraph}
\begin{prop}\label{Prop: Lambda_I -left adjoint}
	The functor {\normalfont $\Lambda_{ \mathfrak{a}}(-): \text{C}(A)_{\mathfrak{a}\text{-cor}}\rightarrow \text{C}(A)_{\mathfrak{a}\text{-red}}$} preserves  all colimits (in particular, coproducts, direct limits, cokernels, pushouts, initial objects and coequilizers). 
\end{prop}

\begin{defn}
	{\normalfont Let $\mathfrak{a}$ be an ideal of $A$ and $M\in\text{C}(A)$.}
	\begin{itemize}
		\item[$(1)$] {\normalfont $M$ is }$\mathfrak{a}$-adically complete {\normalfont if each $M^{n}$ is an $\mathfrak{a}$-adically complete $A$-module \cite[Definition 3.11]{L. Pol-Homotopy theory}}.
		\item [$(2)$] {\normalfont $M$ is} $\mathfrak{a}$-torsion {\normalfont if each $M^{n}$ is an $\mathfrak{a}$-torsion $A$-module}.
	\end{itemize}
\end{defn}
\begin{paragraph}\noi  $\text{C}(A)_{\mathfrak{a}\text{-tor}}$ and $\text{C}(A)_{\mathfrak{a}\text{-com}}$ denote the full subcategory  of $\text{C}(A)$ consisting of $\mathfrak{a}$-torsion and $\mathfrak{a}$-adically complete complexes respectively.
\end{paragraph}
\begin{thm}[The version of MGM equivalence in $\text{C}(A)$]\label{MGM in C(A)} Let $\mathfrak{ a}$ be an ideal of $A$. We have the equalities:
	{\normalfont $$\text{C}(A)_{\mathfrak{a}\text{-tor}}\cap \text{C}(A)_{\mathfrak{a}\text{-red}}= \text{C}(A)_{\mathfrak{a}\text{-com}}\cap \text{C}(A)_{\mathfrak{a}\text{-cor}}=\{M\in\text{C}(A): \mathfrak{ a}M\cong 0\}.$$}
	\end{thm}
\begin{prf}
Let $\mathscr{A}:= \text{C}(A)_{\mathfrak{a}\text{-tor}}\cap \text{C}(A)_{\mathfrak{a}\text{-red}}$ and $\mathscr{B}:= \text{C}(A)_{\mathfrak{a}\text{-com}}\cap \text{C}(A)_{\mathfrak{a}\text{-cor}}$. Let $M\in\mathscr{A}$. Then $\Gamma_{ \mathfrak{a}}(M)\cong M$ and $\mathfrak{ a}\Gamma_{ \mathfrak{a}}(M)\cong 0$ which forces $\mathfrak{ a}M\cong 0$. Conversely, if $M\in\text{C}(A)$ such that $\mathfrak{ a}M\cong 0$, then $\mathfrak{ a}\Gamma_{ \mathfrak{a}}(M)\cong 0$, since $\mathfrak{ a}\Gamma_{ \mathfrak{a}}(M)\subseteq \mathfrak{ a}M$. Moreover, since every element of $M$ is annihilated by $\mathfrak{ a}$, we have $\Gamma_{ \mathfrak{a}}(M)\cong M$. So, $M\in\mathscr{A}$. Thus, $\text{C}(A)_{\mathfrak{ a\text{-red}}}\cap \text{C}(A)_{\mathfrak{ a\text{-tor}}}=\{M\in\text{C}(A): \mathfrak{ a}M\cong 0\}$. Now, let $M\in\mathscr{B}$. We have $M\cong \Lambda_{ \mathfrak{a}}(M)$ and $\mathfrak{ a}\Lambda_{ \mathfrak{a}}(M)\cong 0$. So, $\mathfrak{ a}M\cong 0$. On the other hand, let $M$ be a complex such that $\mathfrak{ a}M\cong 0$. We have $\Lambda_{ \mathfrak{a}}(M)\cong \underset{k}{\varprojlim}M/\mathfrak{ a}^{k}M\cong M$, since $\mathfrak{a}M\supseteq \mathfrak{ a}^{k}M$ for all $k\ge 1$. It follows that $\Lambda_{ \mathfrak{a}}(M)\cong M$, and $\mathfrak{ a}\Lambda_{ \mathfrak{a}}(M)\cong 0$. Thus, $M\in\mathscr{B}$. Therefore, $\text{C}(A)_{\mathfrak{ a\text{-com}}}\cap\text{C}(A)_{\mathfrak{ a\text{-cor}}}=\{M\in\text{C}(A): \mathfrak{ a}M\cong 0\}$.

\end{prf}
\begin{paragraph}\noi
	It is well known that the category of $\mathfrak{ a}$-adically complete modules is not abelian in general, see for instance \cite[Remark 3.3]{porta2014homology}.
\end{paragraph}
\begin{cor} \label{cor- c(A)_com is abelian}Let $\mathfrak{ a}$ be an idempotent ideal of $A$. The full subcategory  {\normalfont $\text{C}(A)_{\mathfrak{ a\text{-com}}}$} of {\normalfont $\text{C}(A)$} is abelian.
	\end{cor}
\begin{prf}
First note that $\text{C}(A)_{\mathfrak{ a\text{-tor}}}$ is abelian, since the category chain complexes of any abelian category is abelian and the category of $\mathfrak{ a}$-torsion modules is abelian. It is enough to show that $\text{C}(A)_{\mathfrak{ a\text{-com}}}=\text{C}(A)_{\mathfrak{ a\text{-tor}}}$. By Remark \ref{Rem: red=cor=any complex}, $\text{C}(A)_{\mathfrak{ a\text{-red}}}=\text{C}(A)_{\mathfrak{ a\text{-cor}}}=\text{C}(A)$. It follows that $\text{C}(A)_{\mathfrak{ a\text{-tor}}}=\text{C}(A)_{\mathfrak{ a\text{-com}}}$, since  by Theorem \ref{MGM in C(A)}, $\text{C}(A)_{\mathfrak{ a\text{-tor}}}\cap \text{C}(A)_{\mathfrak{ a\text{-red}}}= \text{C}(A)_{\mathfrak{ a\text{-tor}}}$, and $\text{C}(A)_{\mathfrak{ a\text{-com}}}\cap \text{C}(A)_{\mathfrak{ a\text{-cor}}}= \text{C}(A)_{\mathfrak{ a\text{-com}}}$. We conclude that $\text{C}(A)_{\mathfrak{ a\text{-com}}}$ is abelian, since $\text{C}(A)_{\mathfrak{ a\text{-tor}}}$ is an abelian category. 
\end{prf}
\section{Results in $\text{D}(A)$}
\begin{paragraph}\noi Let $\mathfrak{ a}$ be an ideal of a ring $A$, and {\normalfont $\text{D}(A)$} denote the derived category of $\text{M}(A)$. In this  section, we prove versions of the GM duality and the MGM equivalence in {\normalfont $\text{D}(A)$} by utilising $\mathfrak{ a}$-reduced and $\mathfrak{ a}$-coreduced complexes. For any complex $M\in\text{C}(A)$, the maps $i:\text{Hom}_{A}(A/\mathfrak{ a},M)\rightarrow \Gamma_{ \mathfrak{a}}(M)$, and $\pi:\Lambda_{ \mathfrak{a}}(M)\rightarrow A/\mathfrak{ a}\otimes_{A}M$ are well-known, see for instance \cite[2.1.3]{Peter-Schenzez}.
\end{paragraph}
\begin{paragraph}\noi Whereas the functor $\Gamma_{ \mathfrak{a}}$ on $\text{C}(A)$ is left exact and therefore has nice homological properties, its dual $\Lambda_{ \mathfrak{a}}$ is not right exact in general, see \cite[Section 1]{Yekutieli-flatness}. To allow for nice homological properties on the functor $\Lambda_{\mathfrak{ a}}$, we impose a condition that makes it to be exact. For a Noetherian ring $A$, $\Lambda_{ \mathfrak{a}}$ is exact on the category of finitely generated $A$-modules,  $\text{M}_\text{f}(A)$ \cite[Section 1]{porta2014homology}, and hence on $\text{C}_\text{f}(A)$, the category of chain of complexes of finitely generated $A$-modules. So, whenever we pass to the local homology modules under these conditions, we will get $\text{L}_{i}\Lambda_{\mathfrak{a}}(M)=0$ for all $i\ne 0$ and $ \text{L}_{0}\Lambda_{\mathfrak{a}}(M)\cong \Lambda_{ \mathfrak{a}}(M)$.
\end{paragraph}
\begin{paragraph}\noi A complex $M\in\text{C}(A)$ is called \textit{acyclic} if $\text{H}^{i}(M)=0$ for all $i\in\Z$. We call a complex $P \in \text{C}(A)$ \textit{$K$-projective} (resp. \textit{$K$-flat}) if for any acyclic complex $M \in \text{C}(A)$ 
	the complex $\text{Hom}_{A} (P, M)$ (resp. $P\otimes_{A}M$) is also acyclic. A complex $ J\in \text{C}(A)$ is called \textit {$K$-injective}
	if for any acyclic complex $M\in \text{C}(A)$ the complex $\text{Hom}_{A}(M, J)$ is also acyclic \cite[Section 2]{porta2014homology}.
%	definitions were introduced in [22]; in [15, Section 3] it is shown that “K-projective”
%	is the same as “having property (P)”, and “K-injective” is the same as “having
%	property (I)”.
	\end{paragraph}

\begin{prop}
	Let {\normalfont $M\in\text{D}(A)$}. The following maps exist in {\normalfont $\text{D}(A)$}.
	{\normalfont \begin{enumerate}
			\item[$(1)$] $\alpha: \text{RHom}_{A}(A/\mathfrak{ a},M)\rightarrow \text{R}\Gamma_{ \mathfrak{a}}(M)$.
			\item[$(2)$] $\beta:\text{L}\Lambda_{ \mathfrak{a}}(M)\rightarrow A/\mathfrak{ a}\otimes_{A}^{\text{L}}M$.
	\end{enumerate}}
\end{prop}
\begin{prf}
	To prove part $(1)$, let $M\rightarrow J$ be a $K$-injective resolution of {\normalfont $M\in\text{D}(A)$}. We have Figure \ref{Figure 1} with $\text{q}_{1},\text{q}_{2}$ quasi-isomorphisms and $i$ an inclusion map.
%\begin{figure}[ht]
%	\centering
%	\begin{minipage}{2in} \begin{tikzcd}
%				\text{RHom}_{A}(A/\mathfrak{ a},M) \arrow{d}{\text{q}_{1}} \arrow[r, dashrightarrow, "\alpha"]
%				& \text{R}\Gamma_{ \mathfrak{a}}(M) \arrow{d}{\text{q}_{2}} \\
%				\text{Hom}_{A}(A/\mathfrak{ a},J) \arrow[r, hook, "i"]
%				& \Gamma_{ \mathfrak{a}}(J)
%			\end{tikzcd}
%		\end{minipage}
%\caption{}\label{Fig 1}
%	\end{figure}
\begin{figure}[ht]
	\begin{center}			 
		%\begin{minipage}{2in}
		\begin{tikzcd} 
			\text{RHom}_{A}(A/\mathfrak{ a},M) \arrow{d}{\text{q}_{1}} \arrow[r, dashrightarrow, "\alpha"]
		& \text{R}\Gamma_{ \mathfrak{a}}(M) \arrow{d}{\text{q}_{2}} \\
		\text{Hom}_{A}(A/\mathfrak{ a},J) \arrow[r, hook, "i"]
		& \Gamma_{ \mathfrak{a}}(J)
		\end{tikzcd}
		%	\end{minipage}
	\caption{}\label{Figure 1}
\end{center}
\end{figure}
	So, take $\alpha:= \text{q}_{2}^{-1}\circ i\circ {q}_{1}$. For part $(2)$, first note that any $K$-flat resolution $F\rightarrow M$ of $M$ induces Figure \ref{Figure 2}
%	\begin{figure}[ht]
%		\centering
%		\begin{minipage}{2in}
%		\begin{tikzcd}
%			\text{L}\Lambda_{ \mathfrak{a}}(M) \arrow{d}{\text{q}_{3}} \arrow[r, dashrightarrow, "\beta"]
%			& A/\mathfrak{ a}\otimes_{A}^{\text{L}}M   \arrow{d}{\text{q}_{4}} \\ \Lambda_{ \mathfrak{a}}(F) 
%			\arrow[r, "\pi"]
%			& A/\mathfrak{ a}\otimes_{A}F
%		\end{tikzcd}
%	\end{minipage}
%\centering \caption{}\label{Figure 2}
%\end{figure}
	\begin{figure}[ht]
	\begin{center}			 
		%\begin{minipage}{2in}
		\begin{tikzcd} 
			\text{L}\Lambda_{ \mathfrak{a}}(M) \arrow{d}{\text{q}_{3}} \arrow[r, dashrightarrow, "\beta"]
			& A/\mathfrak{ a}\otimes_{A}^{\text{L}}M   \arrow{d}{\text{q}_{4}} \\ \Lambda_{ \mathfrak{a}}(F) 
			\arrow[r, "\pi"]
			& A/\mathfrak{ a}\otimes_{A}F	
		\end{tikzcd}
		%	\end{minipage}
	\caption{}\label{Figure 2}
\end{center}
\end{figure}
	where $\text{q}_{3},\text{q}_{4}$ are quasi-isomorphisms and $\pi$ is the canonical surjection. Choose $\beta:= \text{q}_{4}^{-1}\circ \pi\circ \text{q}_{3}$.
	\newline
\begin{paragraph}\noi Recall that for $\star\in\{+,-,\text{b}\}$ we have the full subcategories $\text{C}^{\star}(A)$ (resp. $\text{D}^{\star}(A)$) of $\text{C}(A)$ (resp. $\text{D}(A)$) consisting of bounded below, bounded above and bounded complexes in $\text{C}(A)$ (resp. $\text{D}(A)$). We will write $M\overset{\sim}\rightarrow N$ to denote a quasi-isomorphism between complexes $M$ and $N$.
\end{paragraph}	
\end{prf}
\begin{defn}\label{Defn- red and cor in D(A)}
	{\normalfont A complex $M\in\text{D}^{}(A)$ is} $\mathfrak{ a}$-reduced  {\normalfont if the map $\alpha$ is an isomorphism, i.e., if $\text{RHom}_{A}(A/\mathfrak{ a}, M)\cong\text{R}\Gamma_{ \mathfrak{a}}(M)$}. {\normalfont Dually, a complex $M\in\text{D}^{}(A)$ is}  $\mathfrak{ a}$-coreduced {\normalfont if  $\beta$ is an isomorphism, i.e., if $\text{L}\Lambda_{ \mathfrak{a}}(M)\cong A/\mathfrak{ a}\otimes_{A}^{\text{L}}M$.}
\end{defn}
\begin{paragraph}\noi
	We denote by {\normalfont $\text{D}^{+}(A)_{\mathfrak{ a\text{-red}}}$} (resp. {\normalfont $\text{D}^{-}(A)_{\mathfrak{ a\text{-cor}}}$}) the full subcategories of {\normalfont $\text{D}^{+}(A)$ (resp. $\text{D}^{-}(A)$)} consisting of the bounded below $\mathfrak{ a}$-reduced complexes (resp. the bounded above $\mathfrak{ a}$-coreduced complexes).
\end{paragraph}
\begin{prop}\label{Prop: characterisation of red and cor in D(A)}
	Let $A$ be a ring and {\normalfont $\mathfrak{ a}$} be an ideal of $A$. 
	%The complex  {\normalfont $M\in\text{D}^{+}(A)_{\mathfrak{ a\text{-red}}}$} if and only if {\normalfont $M\in\text{C}^{+}(A)_{\mathfrak{ a\text{-red}}}$}
	
	%	The following conditions are equivalent:
	\begin{enumerate}
		
		\item[$(1)$] The complex  {\normalfont $M\in\text{D}^{+}(A)_{\mathfrak{ a\text{-red}}}$} if and only if  {\normalfont $M\in\text{C}^{+}(A)_{\mathfrak{ a}\text{-red}}$}.
		\item [$(2)$] If $A$ is Noetherian and $M$ is a complex of finitely generated $A$-modules, then the complex  {\normalfont $M\in\text{D}^{-}_\text{}(A)_{\mathfrak{ a\text{-cor}}}$} if and only if  {\normalfont $M\in\text{C}^{-}_\text{}(A)_{\mathfrak{ a}\text{-cor}}$}.
	\end{enumerate}
	\begin{prf}
		{\normalfont	\begin{enumerate}
				\item [$(1)$] Let {\normalfont $M\in\text{D}^{+}(A)$}, and  $M\overset{\sim}\rightarrow J$ be its injective resolution. Then
				{\normalfont 	$M\in\text{D}^{+}(A)_{\mathfrak{ a\text{-red}}}\Leftrightarrow \text{RHom}_{A}(A/\mathfrak{ a},M)\cong \text{R}\Gamma_{ \mathfrak{a}}(M)\Leftrightarrow \text{Hom}_{A}(A/\mathfrak{ a}, J)\overset{\sim}{\rightarrow}\Gamma_{ \mathfrak{a}}(J)\Leftrightarrow \text{H}^{i}(\text{Hom}_{A}(A/\mathfrak{ a}, J))\cong \text{H}^{i}(\Gamma_{ \mathfrak{a}}(J))$} for all $i\in\Z_{\ge 0}$. Let $i\ne 0$. {\normalfont $\text{H}^{i}(\text{Hom}_{A}(A/\mathfrak{ a}, J))=\text{Ext}^{i}_{A}(A/\mathfrak{ a},J)=0\cong \underset{k}{\varinjlim}\text{Ext}^{i}_{A}(A/\mathfrak{ a}^{k},J)=:\text{H}^{i}(\Gamma_{ \mathfrak{a}}(J))$.} 
				For $i=0$, {\normalfont $$\text{H}^{i}(\text{Hom}_{A}(A/\mathfrak{ a}, J))\cong \text{H}^{i}(\Gamma_{ \mathfrak{a}}(J))\Leftrightarrow \text{Hom}_{A}(A/\mathfrak{ a}, M)\cong \Gamma_{ \mathfrak{a}}(M)\Leftrightarrow M\in\text{C}^{+}(A)_{\mathfrak{ a\text{-red}}}.$$}
				\item[$(2)$]  Let $M\in\text{D}^{-}(A)$, and $P\overset{\sim}\rightarrow M$ be its projective resolution. Then
				$M\in\text{D}^{-}_\text{}(A)_{\mathfrak{ a\text{-cor}}}\Leftrightarrow \text{L}\Lambda_{ \mathfrak{a}}(M) \cong A/\mathfrak{ a}\otimes_{A}^\text{L}M\Leftrightarrow \Lambda_{ \mathfrak{a}}(P)\overset{\sim}{\rightarrow}A/\mathfrak{ a}\otimes_{A} P\Leftrightarrow \text{H}_{i}(\Lambda_{ \mathfrak{a}}(P))\cong \text{H}_{i}(A/\mathfrak{ a}\otimes_{A} P)$ for all $i\in\Z_{\leq 0}$. Let $i\ne 0$. $\text{H}_{i}(\Lambda_{ \mathfrak{a}}(P))=0\cong \text{Tor}^{A}_{i}(A/\mathfrak{ a},P)=:\text{H}_{i}(A/\mathfrak{ a}\otimes_{A}P)$. If $i=0$, then  $\text{H}_{i}(\Lambda_{ \mathfrak{a}}(P))\cong \text{H}_{i}(A/\mathfrak{ a}\otimes_{A}P)\Leftrightarrow \Lambda_{ \mathfrak{a}}(M)\cong A/\mathfrak{ a}\otimes_{A}M\Leftrightarrow M\in\text{C}^{-}(A)_{\mathfrak{ a\text{-cor}}}$.
		\end{enumerate} }
	\end{prf}
\end{prop}
\begin{paragraph}\noi By an analogous proof, we get similar results for bounded complexes. That is, $(1)$ the complex $M\in\text{D}^\text{b}(A)_{\mathfrak{ a\text{-red}}}$ if and only if $M\in\text{C}^\text{b} (A)_{\mathfrak{ a\text{-red}}}$, and $(2)$  $M\in\text{D}_\text{}^\text{b}(A)_{\mathfrak{ a\text{-cor}}}$ if and only if $M\in\text{C}_\text{}^\text{b}(A)_{\mathfrak{ a\text{-cor}}}$, provided $A$ is Noetherian and $M$ is a complex of finitely generated $A$-modules.
\end{paragraph}

\begin{exam}\label{Exam: I-red are not closed under q.i} Let $A:=\Z$ and $\mathfrak{ a}:= 2\Z$. In Figure \ref{Figure 3}, {\normalfont $M\in\text{C}(A)_{\mathfrak{ a\text{-red}}}$} and {\normalfont $N\notin \text{C}(A)_{\mathfrak{ a\text{-red}}}$}. However, $f:M\rightarrow N$ is quasi-isomorphism. This shows that, in general, quasi-isomorphisms do not preserve $\mathfrak{ a}$-reduced complexes.

	\begin{figure}[ht]
		\begin{center}			 
		%\begin{minipage}{2in}
		\begin{tikzcd} M:= \cdots\arrow{r}{} & 0\arrow{d}\arrow{r}{}& \Z\arrow{d}\arrow{r}{4}&\Z \arrow{d}\arrow{r}& 0\cdots\\
			N:=	\cdots\arrow{r}{} & 0 \arrow{r}{} & 0 \arrow{r}{}& \Z/4\Z \arrow{r}{}& 0 \cdots	
		\end{tikzcd}
%	\end{minipage}
 \caption{}\label{Figure 3}
\end{center}
\end{figure}	
\end{exam}
 \begin{paragraph}\noi
	Following \cite{vyas2018weak}, we call a complex $M$ \textit{derived $\mathfrak{ a}$-torsion} (resp. \textit{derived $\mathfrak{ a}$-adically complete}) if $\text{R}\Gamma_{ \mathfrak{a}}(M)\cong M$ (resp. $M\cong \text{L}\Lambda_{ \mathfrak{a}}(M)$). Denote by $\text{D}^{+}(A)_{\mathfrak{a}\text{-tor}}$ the full subcategory of $\text{D}(A)$ consisting of bounded below derived $\mathfrak{ a}$-torsion complexes and by $\text{D}^{-}(A)_{\mathfrak{ a}\text{-com}}$ the full subcategory of $\text{D}(A)$ consisting of bounded above derived $\mathfrak{ a}$-adically complete complexes.
	\end{paragraph}
%\begin{paragraph}\noi It is well known that the category of $\mathfrak{ a}$-torsion modules is abelian, and on the other hand, the category of $\mathfrak{ a}$-adically complete modules is not abelian \cite[Remark 3.3]{porta2014homology}. 
%	\end{paragraph}
\begin{paragraph}\noi In \cite[Example 1]{Yekutieli-separated modules}, Yekutieli presented an example of a module which is not $\mathfrak{ a}$-adically complete but derived $\mathfrak{ a}$-adically complete. He, however, proved that a separated derived $\mathfrak{ a}$-adically complete module is $\mathfrak{ a}$-adically complete. In Proposition \ref{Prop: M in derived tor and com iff M in C(A)} $(2)$, we establish that a bounded above complex of finitely generated $A$-modules over a Noetherian ring  is derived $\mathfrak{ a}$-adically complete if and only if it is $\mathfrak{ a}$-adically complete.
	\end{paragraph}
\begin{prop}\label{Prop: M in derived tor and com iff M in C(A)}
	Let $\mathfrak{ a}$ be an ideal of $A$. 
	\begin{enumerate}
		\item [$(1)$] The complex {\normalfont $M\in\text{D}^{+}(A)_{\mathfrak{ a\text{-tor}}}$} if and only if {\normalfont $M\in\text{C}^{+}(A)_{\mathfrak{ a\text{-tor}}}$}.
		\item[$(2)$]If $A$ is Noetherian and $M$ is a complex of finitely generated $A$-modules, then {\normalfont $M\in\text{D}^{-}_\text{}(A)_{\mathfrak{ a\text{-com}}}$} if and only if {\normalfont $M\in\text{C}^{-}(A)_{\mathfrak{ a\text{-com}}}$}.
	\end{enumerate} 
\end{prop}
\begin{prf}
	
	\begin{enumerate}
		\item[$(1)$]  	Let $M\in\text{D}^{+}(A)$.
		$M\in\text{D}^{+}(A)_{\mathfrak{ a\text{-tor}}}\Leftrightarrow \text{R}\Gamma_{ \mathfrak{a}}(M)\cong M \Leftrightarrow  \Gamma_{ \mathfrak{a}}(J)\overset{\sim}{\rightarrow} J\Leftrightarrow \text{H}^{i}(\Gamma_{ \mathfrak{a}}(J))\cong \text{H}^{i}(J)$ for any $i\in\Z_{\ge 0}$, and any injective resolution $M\overset{\sim}\rightarrow J$ of $M$. Let $i\ne 0$. $\text{H}^{i}(\Gamma_{ \mathfrak{a}}(J))=\underset{k}{\varinjlim}\text{Ext}^{i}_{A}(A/\mathfrak{ a}^{k},J)=0\cong \text{H}^{i}(J)$. If $i=0$, then  $\text{H}^{i}(\Gamma_{ \mathfrak{a}}(J))\cong \text{H}^{i}(J)\Leftrightarrow \Gamma_{ \mathfrak{a}}(M)\cong M \Leftrightarrow M\in\text{C}^{+}(A)_{\mathfrak{ a\text{-tor}}}$.
		%		\begin{enumerate}
			%			\item [$(a)$] let $i\ne 0$. $\text{H}^{i}(\Gamma_{ \mathfrak{a}}(J))=\underset{k}{\varinjlim}\text{Ext}^{i}_{A}(A/\mathfrak{ a}^{k},J)=0=\text{H}^{i}(J)$.
			%			\item[$(b)$] if $i=0$, then  $\text{H}^{i}(\Gamma_{ \mathfrak{a}}(J))\cong \text{H}^{i}(J)\Leftrightarrow \Gamma_{ \mathfrak{a}}(M)\cong M \Leftrightarrow M\in\text{C}(A)_{\mathfrak{ a\text{-tor}}}$.
			%		\end{enumerate}
		\item[$(2)$] Let $M\in\text{D}^{-}(A)$ and $P\overset{\sim}\rightarrow M$ be its projective resolution. Then
		$M\in\text{D}^{-}_{\text{}}(A)_{\mathfrak{ a\text{-com}}}\Leftrightarrow M\cong \text{L}\Lambda_{ \mathfrak{a}}(M) \Leftrightarrow P\overset{\sim}{\rightarrow}\Lambda_{ \mathfrak{a}}(P)\Leftrightarrow \text{H}_{i}(P)\cong \text{H}_{i}(\Lambda_{ \mathfrak{a}}(P))$ for any $i\in\Z_{\le0}$. Let $i\ne0$. $\text{H}_{i}(P)=0=\text{H}_{i}(\Lambda_{ \mathfrak{a}}(P))$. If $i=0$, then $\text{H}_{i}(P)\cong M\cong \Lambda_{ \mathfrak{a}}(M)=: \text{H}_{i}(\Lambda_{ \mathfrak{a}}(P))\Leftrightarrow M\in\text{C}^{-}(A)_{\mathfrak{ a\text{-com}}}$. 
		%	 \begin{enumerate}
			%	 	\item [$(a)$] let $i\ne0$. $\text{H}_{i}(P)=0=\underset{k}{\varprojlim}\text{Tor}^{A}_{i}(A/\mathfrak{ a}^{k}, P)$.
			%	 	\item[$(b)$] if $i=0$, then $\text{H}_{i}(P)\cong M\cong \Lambda_{ \mathfrak{a}}(M)=: \text{H}_{i}(\Lambda_{ \mathfrak{a}}(P))\Leftrightarrow M\in\text{C}(A)_{\mathfrak{ a\text{-com}}}$. 
			%	 \end{enumerate}
	\end{enumerate}
\end{prf}
\begin{paragraph}\noi
	Proposition \ref{Prop: characterisation of red and cor in D(A)} (resp. Proposition \ref{Prop: M in derived tor and com iff M in C(A)}) above indicate that the notions of $\mathfrak{ a}$-reduced and $\mathfrak{ a}$-coreduced in $\text{D}(A)$ (resp. the notions of $\mathfrak{ a}$-torsion and $\mathfrak{ a}$-adically complete in  $\text{D}(A)$) under some conditions coincide with notions already studied in Section \ref{ Section: C(A) staff} in the setting of $\text{C}(A)$.
\end{paragraph}
\begin{paragraph}\noi 
	%We note that the full subcategory, $\text{C}(A)_{\mathfrak{ a\text{-com}}}$, of $\text{C}(A)$ consisting of complexes of $\mathfrak{ a}$-adically complete modules is not abelian in general, see \cite[Remark 3.3]{porta2014homology}. Valenzuala G. in 
	Whereas the category, $\text{M}(A)_{\mathfrak{ a\text{-tor}}}$, of $\mathfrak{ a}$-torsion modules is abelian, $\text{M}(A)_{\mathfrak{ a\text{-com}}}$, the category of $\mathfrak{ a}$-adically complete modules, is not abelian in general \cite[Remark 3.3]{porta2014homology}, \cite[Example V.1.1.]{G. Valenzuala}. Valenzuala in \cite[page 71]{G. Valenzuala} constructed an abelian subcategory $\text{M}(A)_{\text{L-com}}$, consisting of $L$-complete modules which contains $\text{M}(A)_{\mathfrak{ a\text{-com}}}$, and showed that $\text{M}(A)_{\mathfrak{ a\text{-tor}}}$ is derived equivalent to $\text{M}(A)_{\text{L-com}}$. 
	In Proposition \ref{Prop: abelian subcategories}, we construct an abelian full subcategory of $\text{C}(A)_{\mathfrak{ a\text{-com}}}$ which consists of $\mathfrak{ a}$-adically complete complexes which are $\mathfrak{ a}$-coreduced or equivalently $\mathfrak{ a}$-torsion complexes which are $\mathfrak{ a}$-reduced, c.f. Theorem \ref{MGM in C(A)}. 
	\end{paragraph}
\begin{prop}\label{Prop: abelian subcategories}
	Let $\mathfrak{ a}$ be an ideal of $A$. 
	\begin{enumerate}
		\item [$(1)$] {\normalfont $\mathcal{A}:= \{M\in \text{C}(A): \mathfrak{ a}M\cong 0\}$} is an abelian category.
		\item[$(2)$]   {\normalfont $ \text{D}(\text{C}(A)_{\mathfrak{ a\text{-tor}}}\cap \text{C}(A)_{\mathfrak{ a\text{-red}}})=\text{D}(\text{C}(A)_{\mathfrak{ a\text{-com}}}\cap \text{C}(A)_{\mathfrak{ a\text{-cor}}})$}.
	\end{enumerate}
\end{prop}
\begin{prf}
It is easy to see that $\mathcal{A}$ is closed under taking submodules, quotients, direct sums and direct products. In particular, for any $\varphi\in\text{Hom}_{\mathcal{A}}(M,N)$, $\text{Ker}\varphi$, $\text{Im} \varphi$, $\text{Coim}(\varphi)$ and $\text{Coker}(\varphi)$ exist and belong to $\mathcal{A}$. Moreover, $\text{Im} \varphi$ is isomorphic to $\text{Coim}(\varphi)$ in $\mathcal{A}$. This proves part $(1)$. Proof of part $(2)$ is immediate, since $\text{C}(A)_{\mathfrak{ a\text{-tor}}}\cap\text{C}(A)_{\mathfrak{ a\text{-red}}}$ and $\text{C}(A)_{\mathfrak{ a\text{-com}}}\cap\text{C}(A)_{\mathfrak{ a\text{-cor}}}$  coincide  with the abelian category $\mathcal{A}$ [Theorem \ref{MGM in C(A)}].
\end{prf}
\begin{paragraph}\noi
	Recall that an inverse system of $A$-modules $\{M_{i}\}_{i\in\mathbb{Z}^{+}}$ is called \textit{pro-zero} if for every $i$ there exists $j\ge i$ such that the homomorphism $M_{j}\rightarrow M_{i}$ is zero. 
	Given a sequence $\textbf{a} = (a_{1} ,\dots , a_{n})$ of elements in $A$, and $i\in \mathbb{Z^{+}}$, we let $\textbf{a}^ {i} :=(a^{i}_{1} ,\dots , a^{i}_{n} )$. A collection of the associated Koszul complexes $\{\text{K}(A;\textbf{a}^{ i} )\}_{ i\in\mathbb{Z^{+}}}$ forms an inverse system \cite{Peter-Schenzez,David-Application-I}.
	{\normalfont Let $A$ be a ring. A ﬁnite sequence $\textbf{a}$ in  $A$ is called} weakly proregular {\normalfont if for	every $q < 0$ the inverse system of $A$-modules $ \{\text{H}^{q}(\text{K}(A;\textbf{a}^{ i}))\}_{ i\in\mathbb{Z^{+}}}$ is pro-zero.}
	{\normalfont An ideal $\mathfrak{a}$ in $A$ is called} \textit{weakly proregular} {\normalfont if it is generated by some weakly proregular sequence \textbf{a}}.
\end{paragraph}
\begin{paragraph}\noi For a weakly proregular ideal $\mathfrak{ a}$ of $A$, Porta et'al proved that
	\begin{thm}{\normalfont \cite[Theorem 7.12]{porta2014homology}}\label{Thm: GM duality wiz WPR}
		Let $\mathfrak{a}$ be a weakly proregular ideal of a ring $A$ and {\normalfont $M, N \in \text{D}(A)$}. Then there exists a natural isomorphism in {\normalfont $\text{D}(A)$} given by
		{\normalfont $$\text{RHom}_{A}(\text{R}\Gamma_{ \mathfrak{a}}(M), N ) \cong \text{RHom}_{R} (M, \text{L}\Lambda_{ \mathfrak{a}}(N)).$$}
	\end{thm} 
\end{paragraph}
\begin{paragraph}\noi
	This theorem says that $\text{R}\Gamma_{ \mathfrak{a}}$ is derived left adjoint to $\text{L}\Lambda_{ \mathfrak{a}}$. We give a condition, in terms of $\mathfrak{ a}$-reduced complexes and $\mathfrak{ a}$-coreduced complexes, for ``$\text{L}\Lambda_{ \mathfrak{a}}$ to be  derived left adjoint to $\text{R}\Gamma_{ \mathfrak{a}}$" [Theorem \ref{Thm-version of GM Duality in D(A)}].
 Denote by $\text{D}(A)_{\mathfrak{ a\text{-red}}}$ (resp. $\text{D}(A)_{\mathfrak{ a\text{-cor}}}$) the full subcategory of $\text{D}(A)$ consisting of $\mathfrak{ a}$-reduced complexes and $\mathfrak{ a}$-coreduced complexes. Moreover, let $\text{D}(A)_{\mathfrak{ a\text{-tor}}}$ (resp. $\text{D}(A)_{\mathfrak{ a\text{-com}}}$) denote the full subcategory of $\text{D}(A)$ consisting of derived $\mathfrak{ a}$-torsion (resp. derived $\mathfrak{ a}$-adically complete) complexes. Theorem \ref{Thm-version of GM Duality in D(A)} gives a version of the GM duality in $\text{D}(A)$.
\end{paragraph}
\begin{thm}[The version of GM duality in $\text{D}(A)$]\label{Thm-version of GM Duality in D(A)} Let $\mathfrak{ a}$ be an ideal of a ring $A$, {\normalfont $M\in\text{D}^{}(A)_{\mathfrak{ a\text{-cor}}}$} and {\normalfont $N\in\text{D}^{}(A)_{\mathfrak{ a\text{-red}}}$}. There exists an isomorphism in {\normalfont $\text{D}(A)$} given by  {\normalfont $$ \text{RHom}_{A}(\text{L}\Lambda_{\mathfrak{ a}}(M), N)\cong \text{RHom}_{A}(M,\text{R}\Gamma_{\mathfrak{ a}}(N)).$$}
	
	%{\normalfont \begin{enumerate}
			%		\item [$(1)$] The functor $\text{R}\Gamma_{ \mathfrak{a}}:\text{D}^{}(A)_{\mathfrak{ a\text{-red}}}\rightarrow \text{D}^{}(A)$ is idempotent.
			%		\item[$(2)$] The functor  $\text{L}\Lambda_{ \mathfrak{a}}: \text{D}^{}(A)_{\mathfrak{ a\text{-cor}}}\rightarrow \text{D}^{}(A)$ is idempotent.
			%		\item[$(3)$]  There exists an isomorphism in $\text{D}(A)$ given by  $$ \text{RHom}_{A}(\text{L}\Lambda_{\mathfrak{ a}}(M), N)\cong \text{RHom}_{A}(M,\text{R}\Gamma_{\mathfrak{ a}}(N)).$$
			%	\end{enumerate}
		%}
\end{thm}
\begin{prf}
	By the Hom-Tensor adjointness formula in the derived category, see for instance \cite[4.4.13]{Peter-Schenzez},  $\text{RHom}_{A}(M, \text{RHom}_{A}(A/\mathfrak{ a}, N))\cong \text{RHom}_{A}(A/\mathfrak{ a}\otimes_{A}^\text{L}M,N )$ for any $M,N\in\text{D}(A)$. Since $M\in\text{D}^{}(A)_{\mathfrak{ a\text{-cor}}}$ and $N\in\text{D}^{}(A)_{\mathfrak{ a\text{-red}}}$ by hypothesis, i.e., $\text{L}\Lambda_{ \mathfrak{a}}(M)\cong A/\mathfrak{ a}\otimes_{A}^\text{L}M$ and $\text{R}\Gamma_{ \mathfrak{a}}(N)\cong \text{RHom}_{A}(A/\mathfrak{ a},N)$, it follows that $\text{RHom}_{A}(M,\text{R}\Gamma_{ \mathfrak{a}}(N))\cong \text{RHom}_{A}(\text{L}\Lambda_{ \mathfrak{a}}(M),N)$ for all $M\in\text{D}(A)_{\mathfrak{ a\text{-cor}}}$, and $N\in\text{D}(A)_{\mathfrak{ a\text{-red}}}$.
\end{prf}

%\begin{paragraph}\noi In the context of $\mathfrak{ a}$-reduced complexes and $\mathfrak{ a}$-coreduced complexes in $\text{D}(A)$,  Theorem \ref{Prop: version of GM Duality} and the weak proregularity condition on the ideal $\mathfrak{ a}$ provide Corollary \ref{cor: derived Gam is both derived left and derived right....}. 
%\end{paragraph}
\begin{cor}\label{cor: derived Gam is both derived left and derived right....}
	Let $\mathfrak{ a}$ be a weakly proregular ideal of $A$, {\normalfont $M\in\text{D}(A)_{\mathfrak{ a\text{-cor}}}$} and {\normalfont $N\in\text{D}(A)_{\mathfrak{ a\text{-red}}}$}.
	\begin{enumerate}
		\item [$(1)$] {\normalfont $ \text{RHom}_{A}(\text{L}\Lambda_{\mathfrak{ a}}(M), N)\cong \text{RHom}_{A}(M,\text{R}\Gamma_{\mathfrak{ a}}(N)).$}
		\item [$(2)$] {\normalfont $\text{RHom}_{A}(\text{R}\Gamma_{ \mathfrak{a}}(M), N ) \cong \text{RHom}_{R} (M, \text{L}\Lambda_{ \mathfrak{a}}(N)).$}
	\end{enumerate}
\end{cor}
\begin{prf} Part $(1)$ of the proof is immediate by Theorem \ref{Prop: version of GM Duality}. Since $\mathfrak{ a}$ is weakly proregular, part $(2)$ follows by  Theorem \ref{Thm: GM duality wiz WPR}. 
\end{prf}

\begin{lem}\label{idempotence of RHom and derived tensor}
	Let {\normalfont $M\in\text{D}(A)$}. 
{\normalfont	\begin{enumerate}
		\item [$(1)$]  $\text{RHom}_{A}(A/\mathfrak{ a}, \text{RHom}_{A}(A/\mathfrak{ a},M))\cong \text{RHom}_{A}(A/\mathfrak{a},M)$.
		\item[$(2)$] $A/\mathfrak{ a}\otimes_{A}^\text{L}(A/\mathfrak{ a}\otimes_{A}^\text{L}M)\cong A/\mathfrak{ a}\otimes_{A}^\text{L}M$.
	\end{enumerate}
}
\end{lem}
\begin{prf}
	First note that $A/\mathfrak{ a}\otimes_{A}^\text{L}A/\mathfrak{ a}\cong A/\mathfrak{ a}$. Part $(1)$ follows by the Hom-Tensor adjointness formula in $\text{D}(A)$, see for instance \cite[4.4.13]{Peter-Schenzez}, and by the associativity of the derived tensor products part $(2)$ holds. 
\end{prf}
%\begin{prop}\label{Prop: Iso b/n RHom and RGam and der tensor and total left derived com}
%	 Let {\normalfont $M\in\text{D}(A)$}.
%	\begin{enumerate}
%		\item [$(1)$] {\normalfont $\text{RHom}_{A}(A/\mathfrak{ a},\text{R}\Gamma_{ \mathfrak{a}}(M))\cong \text{R}\Gamma_{ \mathfrak{a}}(M)$} for any {\normalfont $M\in\text{D}(A)_{\mathfrak{ a\text{-red}}}$}.
%        \item[$(2)$] {\normalfont $A/\mathfrak{ a}\otimes_{A}^\text{L}\text{L}\Lambda_{ \mathfrak{a}}(M)\cong \text{L}\Lambda_{ \mathfrak{a}}(M)$} for any {\normalfont $M\in\text{D}(A)_{\mathfrak{ a\text{-cor}}}$}.
%	\end{enumerate}
%	\end{prop}
%\begin{prf}
%The proof is immediate by  Lemma \ref{idempotence of RHom and derived tensor} and Definition \ref{Defn- red and cor in D(A)}.
%\end{prf}
\begin{lem}\label{rem: RHom in tor iff in red +}
	Let $\mathfrak{ a}$ be an ideal of $A$, and {\normalfont $M\in\text{D}(A)$}. Then {\normalfont $\text{RHom}_{A}(A/\mathfrak{ a},M)\in\text{D}(A)_{\mathfrak{ a\text{-tor}}}$} if and only if  {\normalfont $\text{RHom}_{A}(A/\mathfrak{ a},M)\in\text{D}(A)_{\mathfrak{ a\text{-red}}}$}. Dually, {\normalfont $A/\mathfrak{ a}\otimes_{A}^\text{L}M\in \text{D}(A)_{\mathfrak{ a\text{-com}}}$} if and only if {\normalfont $A/\mathfrak{ a}\otimes_{A}^\text{L}M\in\text{D}(A)_{\mathfrak{ a\text{-cor}}}$}.
\end{lem}
\begin{prf}
%First note that for any $M\in\text{D}(A)$, $\text{RHom}_{A}(A/\mathfrak{ a},\text{RHom}_{A}(A/\mathfrak{ a},M))\cong \text{RHom}_{A}(A/\mathfrak{ a},M)$ and $A/\mathfrak{ a}\otimes_{A}^\text{L}(A/\mathfrak{ a}\otimes_{A}^\text{L}M)\cong A/\mathfrak{ a}\otimes_{A}^\text{L}M$, since the isomorphisms $\text{RHom}_{A}(A/\mathfrak{ a},\text{RHom}_{A}(A/\mathfrak{ a},M))\cong \text{RHom}_{A}(A/\mathfrak{ a}\otimes_{A}^\text{L}A/\mathfrak{ a}, M)\cong \text{RHom}_{A}(A/\mathfrak{ a},M)$ and $A/\mathfrak{ a}\otimes_{A}^\text{L}(A/\mathfrak{ a}\otimes_{A}^\text{L}M)\cong (A/\mathfrak{ a}\otimes_{A}^\text{L}A/\mathfrak{ a})\otimes_{A}^\text{L}M\cong A/\mathfrak{ a}\otimes_{A}^\text{L}M$ are evident.
Suppose that $\text{RHom}_{A}(A/\mathfrak{ a},M)\in \text{D}(A)_{\mathfrak{ a\text{-tor}}}$.  $\text{R}\Gamma_{ \mathfrak{a}}(\text{RHom}_{A}(A/\mathfrak{ a},M))\cong \text{RHom}_{A}(A/\mathfrak{ a},M)\cong \text{RHom}_{A}(A/\mathfrak{ a},\text{RHom}_{A}(A/\mathfrak{ a},M))$, see Lemma \ref{idempotence of RHom and derived tensor} $(1)$. It follows that  $\text{RHom}_{A}(A/\mathfrak{ a},M)\in\text{D}(A)_{\mathfrak{ a\text{-red}}}$ [Definition \ref{Defn- red and cor in D(A)}]. If $\text{RHom}_{A}(A/\mathfrak{ a}, M)\in\text{D}(A)_{\mathfrak{ a\text{-red}}}$, then  $\text{R}\Gamma_{ \mathfrak{a}}(\text{RHom}_{A}(A/\mathfrak{ a},M))\cong\text{RHom}_{A}(A/\mathfrak{ a},\text{RHom}_{A}(A/\mathfrak{ a},M))\cong\text{RHom}_{A}(A/\mathfrak{ a},M)$. The proof of the dual statement is similar.
\end{prf}
\begin{paragraph}\noi For a weakly proregular ideal $\mathfrak{ a}$ of a commutative ring $A$,  Porta  et'al proved the MGM equivalence:
\end{paragraph}
\begin{thm}{\normalfont \cite[Theorem 7.11]{porta2014homology}}
	Let $A$ be a commutative ring, and let $\mathfrak{ a}$ be a weakly proregular ideal in $A$.
	\begin{enumerate}
		\item [$(1)$] For any {\normalfont $M\in\text{D}(A)$, $\text{R}\Gamma_{ \mathfrak{a}}(M)\in\text{D}(A)_{\mathfrak{ a\text{-tor}}}$} and {\normalfont $\text{L}\Lambda_{ \mathfrak{a}}(M)\in \text{D}(A)_{\mathfrak{ a\text{-com}}}$}.
		\item[$(2)$] The functor {\normalfont  $\text{R}\Gamma_{ \mathfrak{a}}: \text{D}(A)_{\mathfrak{ a\text{-com}}}\rightarrow \text{D}(A)_{\mathfrak{ a\text{-tor}}}$} is an equivalence, with quasi-inverse {\normalfont $\text{L}\Lambda_{ \mathfrak{a}}$}. 
	\end{enumerate}
\end{thm}
\begin{prop}\label{Prop:composition} Let $\mathfrak{ a}$ be an ideal of a ring $A$, and {\normalfont $M\in\text{D}^\text{b}(A)$}. 
	\begin{enumerate}
		\item [$(1)$] {\normalfont $\text{RHom}_{A}(A/\mathfrak{ a}, A/\mathfrak{ a}\otimes_{A}^\text{L}M)\cong A/\mathfrak{ a}\otimes_{A}^\text{L}M$} provided $M$ has finite Tor dimension over $A$.
		\item[$(2)$] Let $A$ be Noetherian. {\normalfont $A/\mathfrak{ a}\otimes_{A}^\text{L}\text{RHom}_{A}(A/\mathfrak{ a},M)\cong \text{RHom}_{A}(A/\mathfrak{ a},M)$} provided the complexes $A/\mathfrak{ a}$ and $M$ have finite injective dimension over $A$.
	\end{enumerate}
\end{prop}
\begin{prf}
	Since $M$ has finite Tor dimension,  part $(1)$ of the proof is a consequence of  \cite[Exercise 10.8.3]{weibel1995introduction}, i.e.,  
	$\text{RHom}_{A}(A/\mathfrak{ a}, A/\mathfrak{ a}\otimes_{A}^\text{L}M)\cong \text{RHom}_{A}(A/\mathfrak{ a},A/\mathfrak{ a})\otimes_{A}^\text{L}M\cong A/\mathfrak{ a} \otimes_{A}^\text{L}M$. To prove part $(2)$, consider an injective resolution $M \overset{\sim}\rightarrow J$, where $J$ is a bounded complex of injectives, see \cite[Exercise 10.7.2]{weibel1995introduction}. Also, let $A/\mathfrak{a}\overset{\sim}\rightarrow J^{'}$ be an injective resolution of the complex $A/\mathfrak{ a}$ with $J^{'}\in \text{C}^{b}(A)$. Since $A/\mathfrak{ a}$ has finitely generated cohomology modules ($\text{H}^{i}(A/\mathfrak{ a})=0 $ for $i\ne 0$ and $\text{H}^{i}(A/\mathfrak{ a})=A/\mathfrak{ a}$ for $i=0$, and $A$ is Noetherian), by \cite[Proposition 11.1.5]{Peter-Schenzez} the evaluation morphism $A/\mathfrak{ a}\otimes_{A}\text{Hom}_{A}(J^{'},J)\rightarrow \text{Hom}_{A}(\text{Hom}_{A}(A/\mathfrak{ a},J^{'}),J)$ is a quasi-isomorphism. It follows that $A/\mathfrak{ a}\otimes_{A}^\text{L}\text{RHom}_{A}(A/\mathfrak{ a},M)\cong \text{RHom}_{A}(\text{RHom}_{A}(A/\mathfrak{ a},A/\mathfrak{ a}),M)\cong \text{RHom}_{A}(A/\mathfrak{ a},M)$.
	% So, $A/\mathfrak{ a}\otimes_{A}^\text{L}\text{RHom}_{A}(A/\mathfrak{ a},M) \cong \text{RHom}_{A}(A/\mathfrak{ a},M)$, since {\normalfont $\text{RHom}_{A}(A/\mathfrak{ a},A/\mathfrak{ a})\cong A/\mathfrak{ a}$}.
\end{prf}
\begin{paragraph}\noi
	% Denote by $\text{C}_{\text{f}}(A)$, the category of complexes of finitely generated $A$-modules.  We note that when $A$ is Noetherian, the functor $\Lambda_{ \mathfrak{a}}(-)$ is exact on $\text{C}_{\text{f}}(A)$, since it is exact on $\text{M}_{\text{f}}(A)$, the category of finitely generated $A$-modules, see for instance \cite[page 38]{porta2014homology}, and that $\Lambda_{ \mathfrak{a}}(-)$ extends naturally to complexes. 
	
	 We denote by $\text{D}^{\text{b}}_{\text{f}}(A)$ the full subcategory of $\text{D}(A)$ consisting of bounded complexes of finitely generated $A$-modules.
	\end{paragraph}

\begin{lem}\label{Lem: 4 MGM equality}
	Let $\mathfrak{ a}$ be an ideal of $A$, and {\normalfont $M\in\text{D}^{\text{b}}(A)$}.
	\begin{enumerate}
		\item [$(1)$] {\normalfont $\text{RHom}_{A}(A/\mathfrak{ a},M)\cong M$ if and only if $M\in\text{D}^\text{b}(A)_{\mathfrak{ a\text{-tor}}}\cap \text{D}^\text{b}(A)_{\mathfrak{ a\text{-red}}}$}.
%		\item[$(2)$] Let $\mathfrak{ a}$ be finitely generated. {\normalfont $A/\mathfrak{ a}\otimes_{A}^\text{L}M\cong M$} if and only if {\normalfont $M\in\text{D}^\text{b}(A)_{\mathfrak{ a\text{-com}}}\cap \text{D}^\text{b}(A)_{\mathfrak{ a\text{-cor}}}$}. 
		
			\item[$(2)$] Let $A$ be Noetherian, and {\normalfont $M\in\text{D}^\text{b}_\text{f}(A)$}. Then {\normalfont $A/\mathfrak{ a}\otimes_{A}^\text{L}M\cong M$} if and only if {\normalfont $M\in\text{D}^\text{b}_\text{f}(A)_{\mathfrak{ a\text{-com}}}\cap \text{D}^\text{b}_\text{f}(A)_{\mathfrak{ a\text{-cor}}}$}.
	\end{enumerate}
\end{lem}
\begin{prf}
	\begin{enumerate}
		
	\item[$(1)$]  Suppose  that $M\overset{\sim}\rightarrow  J$  is an injective resolution of the complex $M$.  By \cite[Theorem 1 of section III.7]{Gelfand and Manin}, $\text{H}^{i}(\Gamma_{ \mathfrak{a}}(\text{Hom}_{A}(A/\mathfrak{ a},J)))\cong \text{H}^{i}\Gamma_{ \mathfrak{a}}(\text{H}^{i}\text{Hom}_{A}(A/\mathfrak{ a},J))\cong \text{H}^{i}\Gamma_{ \mathfrak{a}}(\text{Ext}^{i}_{A}(A/\mathfrak{ a},J))\cong 0$ for $i> 0$, and for $i=0$, we have the isomorphism $\text{H}^{i}(\Gamma_{ \mathfrak{a}}(\text{Hom}_{A}(A/\mathfrak{ a},J)))\cong \Gamma_{ \mathfrak{a}}(\text{Hom}_{A}(A/\mathfrak{ a},M))$, see \cite[2.5.1]{weibel1995introduction}. On the other hand, $\text{H}^{i}\text{Hom}_{A}(A/\mathfrak{ a},J)\cong \text{Ext}^{i}_{A}(A/\mathfrak{ a},J)=0$ for $i>0$, and $\text{H}^{0}\text{Hom}_{A}(A/\mathfrak{ a},J)\cong \text{Ext}^{0}_{A}(A/\mathfrak{ a},J)\cong \text{Hom}_{A}(A/\mathfrak{ a},M)$ which is $\mathfrak{ a}$-torsion.  So, $\Gamma_{ \mathfrak{a}}(\text{Hom}_{A}(A/\mathfrak{ a},J))\overset{}\rightarrow \text{Hom}_{A}(A/\mathfrak{ a},J)$ is a quasi-isomorphism, and thus $\text{RHom}_{A}(A/\mathfrak{ a},M)\in\text{D}^\text{b}(A)_{\mathfrak{ a\text{-tor}}}$.  It follows that $\text{RHom}_{A}(A/\mathfrak{ a},M)\in\text{D}^\text{b}(A)_{\mathfrak{ a\text{-red}}}$ [Lemma \ref{rem: RHom in tor iff in red +}]. Thus, $M\in\text{D}^\text{b}(A)_{\mathfrak{ a\text{-tor}}}\cap \text{D}^\text{b}(A)_{\mathfrak{ a\text{-red}}}$, since, by hypothesis, $\text{RHom}_{A}(A/\mathfrak{ a},M)\cong M$.  If $M\in\text{D}^\text{b}(A)_{\mathfrak{ a\text{-tor}}}\cap \text{D}^\text{b}(A)_{\mathfrak{ a\text{-red}}}$, the isomorphisms $\text{R}\Gamma_{ \mathfrak{a}}(M)\cong M$, and $\text{R}\Gamma_{ \mathfrak{a}}(M)\cong \text{RHom}_{A}(A/\mathfrak{ a},M)$, yield $\text{RHom}_{A}(A/\mathfrak{ a},M)\cong M$.
	\item[$(2)$] Choose a projective resolution $F\overset{\sim}\rightarrow M$ of the complex $M$. Since  $\Lambda_{ \mathfrak{a}}(-)$ is right exact on $\text{C}^\text{b}_\text{f}(A)$, by the dual version of \cite[Theorem 1 of section III.7]{Gelfand and Manin}, $\text{H}_{i}(\Lambda_{ \mathfrak{a}}(A/\mathfrak{ a}\otimes_{A} F))\cong \text{H}_{i}\Lambda_{ \mathfrak{a}}(\text{H}_{i}(A/\mathfrak{ a}\otimes_{A} F)\cong \text{H}_{i}\Lambda_{ \mathfrak{a}}(\text{Tor}^{A}_{i}(A/\mathfrak{ a},F))\cong 0$ for $i> 0$, and for $i=0$, we have an isomorphism $\text{H}_{i}\Lambda_{ \mathfrak{a}}(A/\mathfrak{ a}\otimes_{A}F)\cong \Lambda_{ \mathfrak{a}}(A/\mathfrak{ a}\otimes_{A}M))$, see \cite[2.5.1]{weibel1995introduction}. On the other hand, there are isomorphisms $\text{H}_{i} (A/\mathfrak{ a}\otimes_{A}F)\cong \text{Tor}_{i}^{A}(A/\mathfrak{ a},F)=0$ for $i>0$, and $\text{H}_{0}(A/\mathfrak{ a}\otimes_{A}F)\cong \text{Tor}_{0}^{A}(A/\mathfrak{ a}\otimes_{A}F)\cong A/\mathfrak{ a}\otimes_{A}M$ which is $\mathfrak{ a}$-adically complete, and also finitely generated, since $A/\mathfrak{ a}$ and $M$ are finitely generated by hypothesis.  So, the morphism $\Lambda_{ \mathfrak{a}}(A/\mathfrak{ a}\otimes_{A}F)\overset{}\rightarrow A/\mathfrak{ a}\otimes_{A}F$ is a quasi-isomorphism, and thus $A/\mathfrak{ a}\otimes_{A}^\text{L}M\in\text{D}_\text{f}^\text{b}(A)_{\mathfrak{ a\text{-com}}}$.  It follows that $A/\mathfrak{ a}\otimes_{A}^\text{L}M\in\text{D}_\text{f}^\text{b}(A)_{\mathfrak{ a\text{-cor}}}$ [Lemma \ref{rem: RHom in tor iff in red +}]. So, $M\in\text{D}^\text{b}_\text{f}(A)_{\mathfrak{ a\text{-com}}}\cap \text{D}_\text{f}^\text{b}(A)_{\mathfrak{ a\text{-cor}}}$, since, by hypothesis, $A/\mathfrak{ a}\otimes_{A}^\text{L}M\cong M$. Now, if $M\in\text{D}^\text{b}_\text{f}(A)_{\mathfrak{ a\text{-com}}}\cap \text{D}^\text{b}_\text{f}(A)_{\mathfrak{ a\text{-cor}}}$, then from the isomorphisms $M\cong \text{L}\Lambda_{ \mathfrak{a}}(M)$, and $\text{L}\Lambda_{ \mathfrak{a}}(M)\cong A/\mathfrak{ a}\otimes_{A}^\text{L}M$, it follows that  $A/\mathfrak{ a}\otimes_{A}^\text{L}M\cong M$.
	\end{enumerate}
\end{prf}
\begin{thm}[The version of MGM equivalence in $\text{D}(A)$]\label{Theorem- MGM equality}
	Let $A$ be a Noetherian ring such that each complex in {\normalfont $\text{D}_\text{f}^\text{b}(A)$} has both finite injective dimension and finite Tor dimension. We have the equality:
	{\normalfont \begin{equation*}
			\text{D}_\text{f}^\text{b}(A)_{\mathfrak{a}\text{-com}}\cap \text{D}_\text{f}^\text{b}(A)_{\mathfrak{a}\text{-cor}}= \text{D}_\text{f}^\text{b}(A)_{\mathfrak{a}\text{-tor}}\cap \text{D}_\text{f}^\text{b}(A)_{\mathfrak{a}\text{-red}}.
	\end{equation*}}
\end{thm}
\begin{prf}
	Let {\normalfont $M\in \text{D}_\text{f}^\text{b}(A)_{\mathfrak{a}\text{-com}}\cap \text{D}_\text{f}^\text{b}(A)_{\mathfrak{a}\text{-cor}}$}. Since $M$ has finite Tor dimension, the isomorphism $\text{RHom}_{A}(A/\mathfrak{ a},M)\cong M$ follows from Proposition \ref{Prop:composition} $(1)$, i.e., $\text{RHom}_{A}(A/\mathfrak{a},A/\mathfrak{ a}\otimes_{A}^\text{L}M)\cong A/\mathfrak{ a}\otimes_{A}^\text{L}M$, since by definition,  $M\cong\text{L}\Lambda_{ \mathfrak{a}}(M)$ and $\text{L}\Lambda_{ \mathfrak{a}}(M)\cong A/\mathfrak{ a}\otimes_{A}^\text{L}M$ . The conclusion $M\in\text{D}_\text{f}^\text{b}(A)_{\mathfrak{a}\text{-tor}}\cap \text{D}_\text{f}^\text{b}(A)_{\mathfrak{a}\text{-red}}$ follows by Lemma \ref{Lem: 4 MGM equality} $(1)$. Conversely, let $M\in \text{D}_\text{f}^\text{b}(A)_{\mathfrak{a}\text{-tor}}\cap \text{D}_\text{f}^\text{b}(A)_{\mathfrak{a}\text{-red}}$ have finite injective dimension. Since $\text{R}\Gamma_{ \mathfrak{a}}(M)\cong M$, and $\text{R}\Gamma_{ \mathfrak{a}}(M)\cong \text{RHom}_{A}(A/\mathfrak{ a},M)$  and by Proposition \ref{Prop:composition} $(2)$, i.e., $A/\mathfrak{ a}\otimes_{A}^\text{L}\text{RHom}_{A}(A/\mathfrak{ a},M)\cong \text{RHom}_{A}(A/\mathfrak{ a},M)$, the isomorphism $A/\mathfrak{ a}\otimes_{A}^\text{L}M\cong M$ is immediate, we conclude that $M\in\text{D}_{\text{f}}^\text{b}(A)_{\mathfrak{ a\text{-com}}}\cap \text{D}^\text{b}_{\text{f}}(A)_{\mathfrak{ a\text{-cor}}}$ by Lemma \ref{Lem: 4 MGM equality} $(2)$.
\end{prf}
\begin{cor}\label{MGM equality}
	Let the ideal $\mathfrak{ a}$ of a Noetherian ring $A$ be idempotent such that each complex in {\normalfont $\text{D}_{\text{f}}^\text{b}(A)$} has finite injective dimension and finite Tor dimension. We have {\normalfont	$$\text{D}_{\text{f}}^\text{b}(A)_{\mathfrak{a}\text{-com}}= \text{D}_\text{f}^\text{b}(A)_{\mathfrak{a}\text{-tor}}.$$}
\end{cor}
\begin{prf}
	Since $\mathfrak{ a}=\mathfrak{ a}^{2}$, we get $\text{C}^{\text{b}}_{}(A)=\text{C}_{}^\text{b}(A)_{\mathfrak{ a\text{-cor}}}=\text{C}^{\text{b}}_{}(A)_{\mathfrak{ a\text{-red}}}$. By  virtue of Proposition \ref{Prop: characterisation of red and cor in D(A)} $(2)$ and Proposition \ref{Prop: M in derived tor and com iff M in C(A)} $(2)$,  $\text{D}^{\text{b}}_{\text{f}}(A)=\text{D}_{\text{f}}^\text{b}(A)_{\mathfrak{ a\text{-cor}}}=\text{D}^{\text{b}}_{\text{f}}(A)_{\mathfrak{ a\text{-red}}}$. The conclusion follows from Theorem \ref{Theorem- MGM equality}.
\end{prf}
\begin{paragraph}\noi In conclusion, we have provided conditions, namely; $\mathfrak{ a}$-reduced complexes and $\mathfrak{ a}$-coreduced complexes for which versions of the GM duality and the MGM equivalence hold. Note that these conditions and the condition common in the literature, namely; weakly proregular ideal are independent. None implies the other in general, see \cite[Section 4.1]{David-Application-I}. Finally, this paper just like \cite{reduced w.r.t another-ours,  kyomuhangi2020locally,Annet-David : Generalized reduced, David-Application-I, ssevviiri: App II} has benefited from the fact that the notion of $\mathfrak{ a}$-reduced modules or complexes is categorical, dualisable and amenable to study by use of methods from category theory and homological algebra.
	\end{paragraph}
\section*{Acknowledgment}	
	\begin{paragraph}\noi
	We  acknowledge support from the Eastern Africa Algebra Research Group (EAALG), the EMS-Simons for Africa Program and the International Science Program (ISP). Part of this work was written while the second author was visiting the third author at Makerere University. The second author is grateful for the hospitality.
\end{paragraph}
	
	\end{document}